\documentclass[review]{elsarticle}

\usepackage{lineno,hyperref,eurosym}

\journal{Acta Astronautica}

\usepackage{newtxtext,newtxmath} 
\usepackage{booktabs}
\usepackage{longtable}
\usepackage{geometry}
\usepackage{makecell}
\usepackage{pdflscape} 
\usepackage[table]{xcolor}
\usepackage{colortbl}
\definecolor{LightRed1}{RGB}{255,220,220}
\usepackage{rotating}
\usepackage{amsmath}
\usepackage{nicefrac}
\usepackage{subcaption}
\usepackage{booktabs, makecell, longtable, array, xcolor, rotating}
\usepackage{microtype}
\usepackage[english]{babel}
\usepackage[utf8]{inputenc}
\usepackage{bm}
\usepackage{amsmath}
\usepackage{float}
\usepackage{hyperref}
\usepackage{caption}
\usepackage{graphicx}
\usepackage{multirow}
\usepackage{cleveref}

 \hypersetup{
    colorlinks,
    linkcolor={black!50!},
    citecolor={black},
    urlcolor={black}}

\begin{document}

\begin{frontmatter}

\title{{Low-Energy and Low-Thrust Exploration Tour of Saturnian Moons with Full Lunar Surface Coverage}}





\author[label1]{Chiara Pozzi}
\author[label2]{Mauro Pontani}
\author[label1]{Alessandro Beolchi}
\author[label3]{Hadi Susanto}
\author[label4]{Elena Fantino$^*$}
\address[label1]{Department of Aerospace Engineering, Khalifa University of Science and Technology, P.O. Box 127788, Abu Dhabi, United Arab Emirates}
\address[label2]{Department of Mechanical and Aerospace Engineering, Sapienza Università di Roma, Rome, Italy}
\address[label3]{Department of Mathematics, Khalifa University of Science and Technology, P.O. Box 127788, Abu Dhabi, United Arab Emirates}
\address[label4]{Department of Aerospace Engineering and Polar Research Center, Khalifa University of Science and Technology, P.O. Box 127788, Abu Dhabi, United Arab Emirates, elena.fantino@ku.ac.ae}
\cortext[cor]{Corresponding author: elena.fantino@ku.ac.ae (Elena Fantino)}

\begin{abstract}
This study presents the trajectory design for a mission touring the Inner Large Moons of Saturn—Rhea, Dione, Tethys, Enceladus, and Mimas—engineered to satisfy observational requirements, including full surface coverage, as defined by the scientific community and major space agencies, while ensuring low fuel consumption and compatibility with current power and propulsion technologies (radioisotope thermoelectric generators and Hall effect thrusters). The tour begins at Rhea and concludes at Mimas, employing a trajectory concept that alternates between extended observation phases around each moon and Saturn-centered low-thrust spiral arcs to transition efficiently to the next target in the sequence. The $J_2$-perturbed Circular Restricted Three-Body Problem is the dynamical model adopted for the design of exploration paths around the moons, with halo orbits used as staging points for heteroclinic and homoclinic loops that guarantee prolonged, repeated, and comprehensive surface reconnaissance of the targets (including critical regions such as Enceladus’ poles, where geological activity manifests through intense material plumes). Stable and unstable hyperbolic invariant manifolds of the halo orbits act as departure and arrival gateways for propelled inter-moon transfers. The latter are modeled in an ephemeris-based framework with the inclusion of the relevant gravitational perturbations induced by the moons of the system, the Sun, and the oblateness of Saturn. The dynamical model setup is carried out through a rigorous perturbation analysis to maximize computational efficiency while ensuring a high-fidelity trajectory design. A locally-optimal guidance law is used to minimize propellant consumption. The proposed tour represents an alternative to traditional flyby-based missions, offering comparable total duration yet with a greater fraction of observing time and reduced fuel requirements. It advances previous work by achieving both complete lunar surface coverage and high-fidelity trajectory modeling.
\end{abstract}
\end{frontmatter}


\section{Introduction}
The moons of giant planets are among the most diverse and scientifically compelling objects in the Solar System. Among them, Saturn stands out for its exceptional group of major moons (Mimas, Enceladus, Tethys, Dione, Rhea), known as the Inner Large Moons (ILMs). The Cassini-Huygens mission marked a turning point in the exploration of the Ringed Planet and its complex lunar system \cite{matson2002cassini, spilker2019cassini}. In addition, it has brought to light several scientific questions that remain open. As stated in the Decadal Strategy for Planetary Science and Astrobiology 2023-2032 \cite{NAP26522}, answering them requires in situ exploration of the giant planetary system. This statement aligns with the goals of the ESA Voyage 2050 Program, which prioritizes large-class science missions aimed at exploring the moons of the giant planets, with a particular emphasis on Saturn and Enceladus \cite{ESA_Moons_Giant_Planets_2024}. Ocean worlds in the outer Solar System are particularly captivating targets for scientific investigation, as they are key for understanding the evolution of the Solar System and represent prime candidates in the search for biosignatures. Among the ILMs, Enceladus stands out for its subsurface activity and geyser-like jets at its poles, which provide direct access to materials originating beneath the icy crust \cite{mackenzie2016theo, miles2025endogenic}. For this reason, the ESA Voyage 2050 report \cite{ESA_Moons_Giant_Planets_2024} identifies an Enceladus orbiter mission as one of the most urgent scientific priorities. Other ILMs, such as Mimas and Dione, also hold significant scientific potential. These primarily icy bodies are believed to possibly host a subsurface ocean, but conclusive evidence is still missing \cite{buratti2018cold, dougherty2018review}. Similarly, the origin of the red streaks on Tethys is unknown, underscoring their value as targets for in-depth exploration. The report further emphasizes the need for high-inclination and low-altitude orbits to enable comprehensive gravity and magnetometry measurements, as well as detailed surface topography mapping. This highlights the critical importance of dedicated mission for Enceladus as well as the entire system of ILMs.\\
\indent The first conceptual study of a moon tour (MT) was the Petit Grand Tour of the Galilean moons \cite{koon2002constructing, gomez2003libration,koon2000heteroclinic}, which employed dynamical system methods to design transfers between libration point orbits (LPOs) within the Circular Restricted Three-Body Problems (CR3BPs) formed by Jupiter and individual moons. By coupling distinct CR3BPs, the study demonstrated the feasibility of patching arcs across these systems. Simplified two-body models (2BP) are frequently employed in trajectory design, with the {V-infinity Leveraging Transfer (VILT)} emerging as a particularly effective method for constructing low-fidelity solutions for moon-to-moon transfers. A VILT is a trajectory that connects two consecutive gravity assists at the same moon, incorporating a {V-infinity Leveraging Maneuver (VILM)}. The VILTs are used to design MTs through {resonance hopping}, which patches consecutive orbits in mean motion resonance with a given moon. This approach has been developed by Campagnola et al.~\cite{campagnola2010endgame, campagnola2010fast} to design low-$\Delta V$ orbital insertion maneuvers around the moons of Jupiter and Saturn. Strange, Campagnola, and Russell~\cite{strange2009leveraging} developed a leveraging tour that efficiently enabled capture at Enceladus by combining gravity assists from Titan, Rhea, Dione, and Tethys with VILMs. This approach reduced the propellant consumption to less than half of that required for a direct transfer from Titan, assuming a chemical propulsion system. Resonance hopping can be extended to the CR3BP framework. In fact, different CR3BPs were linked through the Tisserand–Poincaré (T-P) graph, and this led to identifying families of periodic orbits, in the context of a resonant hopping tour. Proposed in~\cite{campagnola2010endgame2} for the Jovian system, the T–P graph method identifies ballistic endgames and links different CR3BP systems to build efficient transfers between moons. Compared to impulsive maneuvers, the T–P graph method significantly reduces propellant consumption by leveraging resonant orbits and high-altitude flybys. Further investigations by Lantoine et al.~\cite{ lantoine2011near, lantoine2011optimization} explored unstable resonant orbits to identify indirect transfers between adjacent moons within the Jovian system. Similar studies focusing on the Saturnian system have also been conducted~\cite{brown2008multi,russell2009cycler,lara2010mission}. More recently, Takubo et al.~\cite{takubo2024automated} completed a heuristically optimized database of VILT transfers between Titan and Enceladus using a full-ephemeris model. 

\begin{figure*}[h]
    \centering
    \includegraphics[width=1\linewidth]{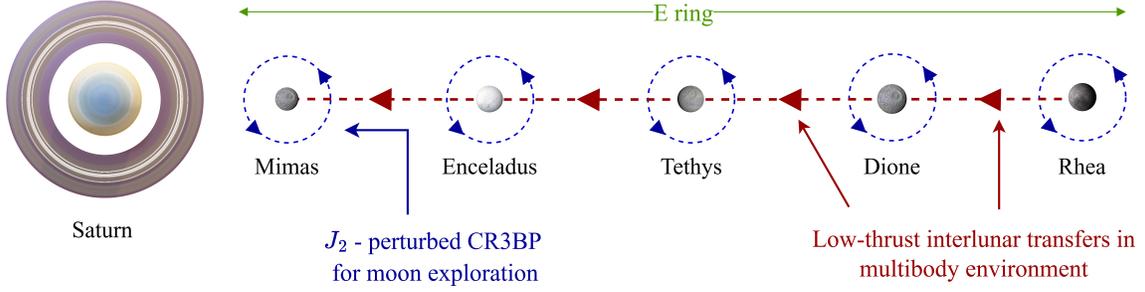}
    \caption{Sketch of the MT strategy, highlighting the dynamical models used in the distinct phases of the tour.}
    \label{fig:sketched}
\end{figure*}

\indent Trajectory design for MTs presents significant dynamical challenges. The Cassini mission relied on a sequence of fast flybys enabled by a chemical-propulsion system on highly elliptical Saturn-centred orbits. During its 13-year orbital tour, approximately 700 kg of propellant were used, corresponding to a total $\Delta V$ of roughly 800 m/s, yet no orbital capture around any ILM was achieved. With chemical propulsion, direct orbital capture at the ILMs requires large propellant consumption, making in-orbit operations unfeasible. An improvement was achieved through the leveraging tour designed in \cite{strange2009leveraging}, which achieved capture at Enceladus via a sequence of VILMs, requiring a total $\Delta V$ of 960~m/s and about 1000~kg of propellant, compared with the 2750~kg needed for a direct transfer from Titan. However, this approach still relies on fast flybys that provide short visibility windows and do not allow extended observational periods. A notable strategy to overcome these challenges was proposed in \cite{fantino2023end}. The method leverages the natural dynamics of the three-body problem and electric propulsion to bridge the energetic gaps that prevent conventional orbital insertion around the Saturnian moons. The approach relies on weak lunar capture through low-energy trajectories, enabling long-duration dwelling, thus allowing investigations for astrobiological research. In certain planetary systems, Hyperbolic Invariant Manifolds (HIMs) associated with LPOs in the CR3BP serve as a powerful tool for designing low-cost MTs. Traditionally, moon-to-moon trajectories are identified through intersections of HIMs using Poincaré sections. However, in the Saturnian system, such intersections do not exist in configuration space due to the near-circular Saturn-centered orbits of the manifolds. For the Saturn case, the MT was studied in two dimensions using a CR3BP-2BP patched method \cite{fantino2019connecting,fantino2023end}. This approach leverages Planar Lyapunov Orbits (PLOs) of the Saturn-moon system as terminal orbits. To connect consecutive moons, the HIMs associated with the PLOs were used. The strategy outlined in \cite{fantino2023end} is based on two-dimensional low-thrust (LT) arcs that link the unstable and stable branches of the invariant manifolds with marginal mass consumption. In particular, the resulting 1.7-year tour, which did not include the largest ILM (Rhea), required a total $\Delta V$ of approximately 3.1 km/s and a propellant usage of about 125 kg. For each moon, exploration relied on low-energy trajectories arising in the CR3BP, specifically two-dimensional homoclinic and heteroclinic connections between PLOs \cite{salazar2021observational}.  Then the science orbits were extended to three-dimensional space by substituting the planar Lyapunov orbits (PLOs) with halo orbits, which allowed the computation of 3D low-energy connections \cite{fantino2020design}. However, the moon-to-moon phases remained confined to the planar approximation. For the exploration orbits, further refinement was achieved by including the effects of the polar flattening of Saturn and Enceladus on the spacecraft dynamics \cite{salazar2021science}. Approximating Saturn and its moons as spherical bodies can introduce inaccuracies, as Saturn is the most oblate planet in the Solar System and many of its moons also deviate from sphericity.
Other works examined the influence of the asphericity of one or both primaries within the CR3BP, showing that it can alter the periodic orbits obtained in the unperturbed framework \cite{bury2020effect}. In \cite{mittal2020analysis, safiya2012oblateness}, the authors analyzed the variations in the geometry and energy of periodic and quasi-periodic orbit families in the planar CR3BP, while \cite{abouelmagd2015effect,abouelmagd2012existence} examined how the oblateness affects the stability of the libration points.
\\
 \begingroup
\indent This work builds on the studies presented in \cite{fantino2023end,salazar2021science}, aiming to overcome their main limitations and advance toward a high-fidelity, fuel-efficient strategy for the exploration of the ILMs of Saturn through a three-dimensional tour. In particular, earlier studies (i) rely on planar, unperturbed CR3BP models to design science orbits, (ii) treat the inter-moon transfers within a 2BP framework rather than an integrated $N$-body model, and (iii) do not provide an end-to-end tour architecture visiting all five ILMs under realistic low-thrust and power constraints.
In this study, the adopted propulsion system is electric, powered by radioisotope thermoelectric generators (RTGs) due to the low solar irradiance in the Saturnian region. This imposes strict constraints on the thrust level; hence, the spacecraft operates with low thrust compatible with attainable electric propulsion capabilities. Motivated by the above gaps, the present work has the following objectives:

\begin{itemize}
    \item extend previous works by transitioning from a planar to a spatial
    formulation while employing a $J_2$-perturbed CR3BP model within the SOI of
    each ILM, leveraging the resulting halo orbits and associated HIMs to build
    low-energy science orbits around the $L_1$ and $L_2$ equilibrium points, with
    the aim of long-duration dwelling and improved coverage of the polar regions;

    \item design three-dimensional moon-to-moon transfers in an $N$-body,
    ephemeris-based dynamical framework, combining low-thrust electric propulsion
    with HIMs, to obtain energetically favorable transfers;

    \item perform a dedicated perturbation analysis to identify and retain only
    those effects that exceed prescribed accuracy thresholds, to reduce the
    computational cost of the inter-moon legs while maintaining the required
    fidelity;

    \item assemble and assess a representative end-to-end tour architecture that
    includes all five ILMs, most notably Rhea, which had not been considered in
    previous studies.
\end{itemize}

\noindent The overall goal of this study is to propose a new exploration strategy of the Saturnian system with minimal propellant usage, by combining weak capture, low-energy trajectories with LT transfer arcs. 
An overview of the proposed MT is shown in Fig.~\ref{fig:sketched}. The trajectory is structured into sequential blocks, from Rhea to Mimas, each consisting of an observation phase around the target moon, followed by an inter-moon powered segment. 

\endgroup


\indent The paper is organized as follows. Section \ref{dynamical} introduces the reference frames and dynamical models. Section \ref{model} presents the analysis of orbit perturbations, identifying and ranking the relevant contributions included in the dynamical model. Section \ref{tourdesign} details the lunar tour design methods, including science orbit construction and inter-moon transfer strategies. Section \ref{resultssec} presents results for the two MT phases. Sections \ref{disc} and \ref{concl} conclude the paper by remarking on the main findings and outlining directions for future work. A preliminary version of this work was presented at the $76^{th}$ International Astronautical Congress \cite{pozzi2025IAC}. The following abbreviations are often used in tables and figures: S = Saturn,
Ti = Titan, Rh = Rhea, Di = Dione, Te = Tethys, En = Enceladus, Mi = Mimas.

\section{Dynamical models} \label{dynamical}
This section outlines the dynamical frameworks employed in the study. The reference frames used throughout the analysis are introduced. The formulation of the $J_2$-perturbed CR3BP, which forms the basis for the design of exploration orbits, is presented next. Finally, the dynamical model adopted for inter-moon transfers is described.

\subsection{Reference frames}\label{refframesection}
This study employs both fixed-axes (denoted by $\mathcal{N}$) and rotating (identified as $\mathcal{R}$) reference frames to describe the spacecraft motion.
The fixed-axes frames are non-rotating with respect to the fixed stars and are therefore treated as inertial for the purposes of this analysis. The Saturn-Centered Inertial (SCI) reference frame ($\mathcal{N}_S ; {O_S; X_S, Y_S, Z_S}$) has its origin at the center of mass $O_S$ of Saturn, and is associated with the right-handed sequence of three mutually orthogonal unit vectors $ \{{\boldsymbol {{c}}_1}^{S}, \: {\boldsymbol {{c}}_2}^{S}, \: {\boldsymbol {{c}}_3}^{S}\}$. Specifically, $\boldsymbol{c}_1^S$ is aligned with the line of nodes defined by the intersection of the orbital plane of Enceladus with the equatorial plane of Saturn at epoch J2000, $\boldsymbol{c}_3^S$ is directed along the rotation axis of Saturn, and $\boldsymbol{c}_2^S$ completes the right-handed triad.\\
\begin{figure*}[t]
    \centering
   \begin{subfigure}[b]
   {0.5\textwidth}
        \centering
        \includegraphics[width=\textwidth]{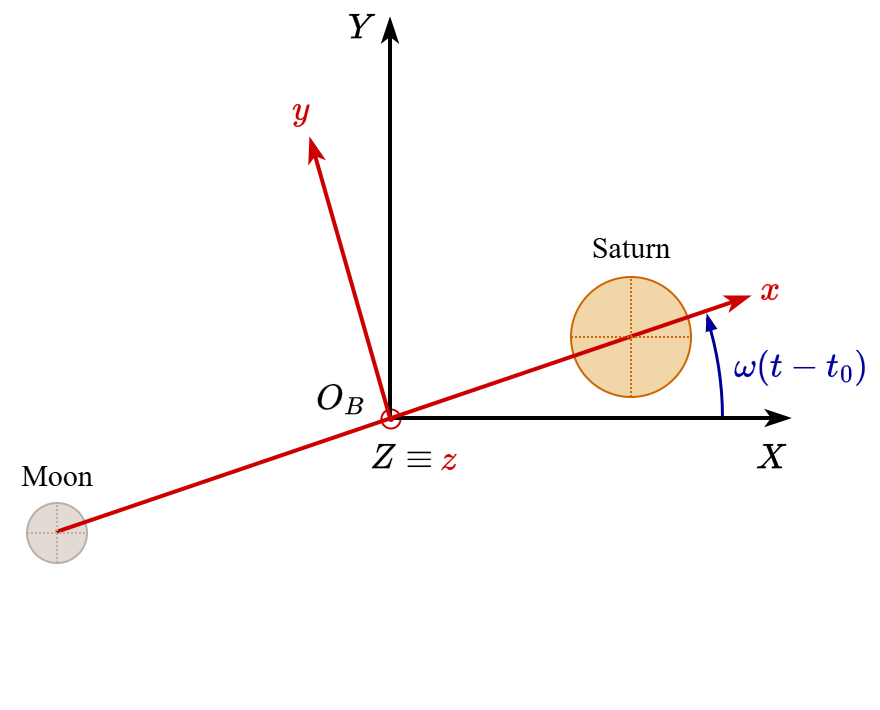} 
        \captionsetup{skip=15pt}
        \caption{Barycentric inertial (black) and synodic (red) reference frames for a Saturn-moon system.}
        \label{SYNframe}
    \end{subfigure}
    \hfill
     \begin{subfigure}[b]{0.45\textwidth}
        \centering
    \includegraphics[width=\textwidth]{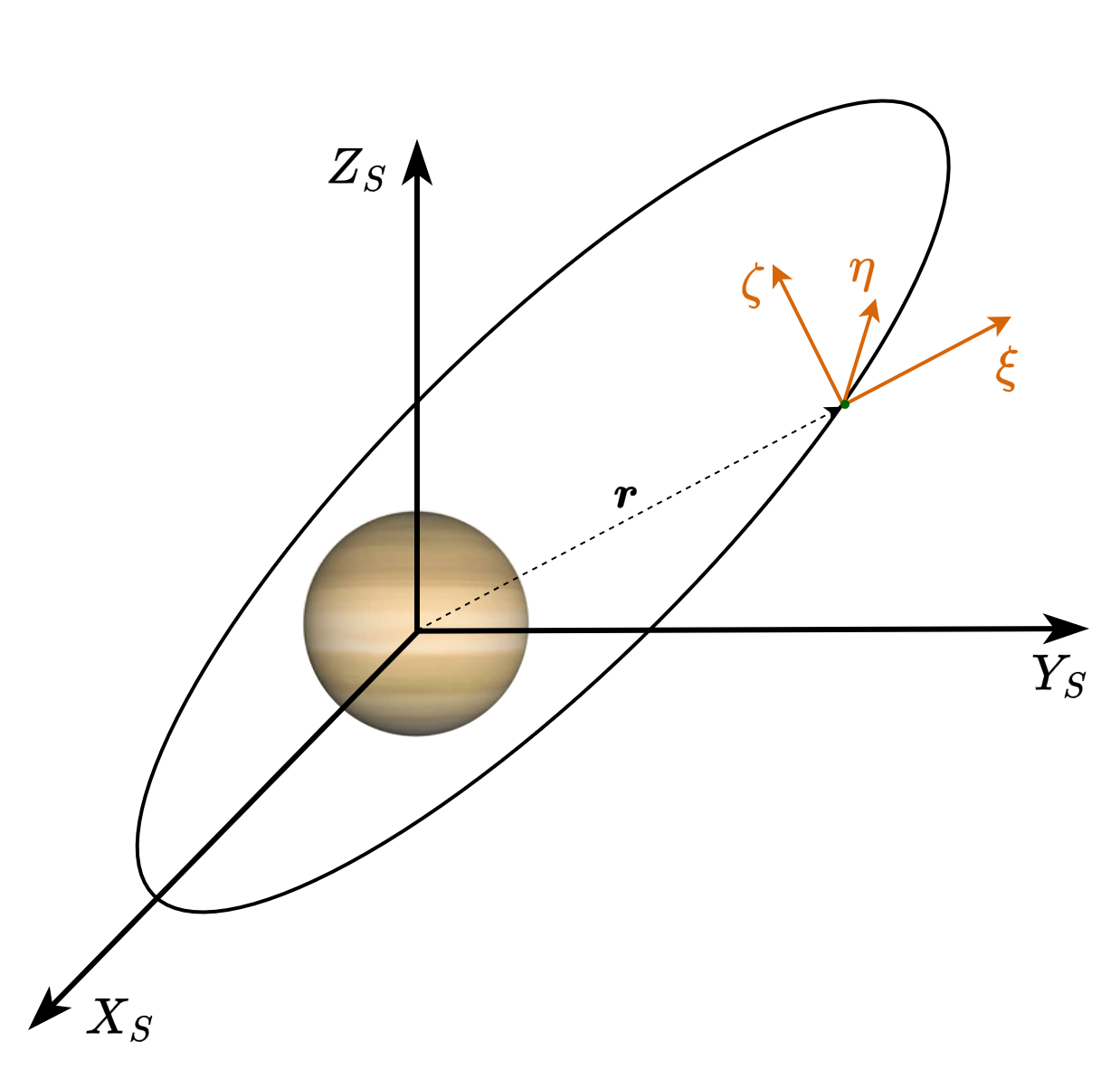} 
        \caption{SCI and spacecraft-centered RSW frame.}
        \label{rswframepic}
    \end{subfigure}

    \caption{{Reference frames used in the dynamical model.}}
    \label{fig:frames}
\end{figure*}
\indent The barycentric inertial frame ($\mathcal{N}_B ; {O_B; X, Y, Z}$) has its origin at the barycenter $O_B$ of a Saturn-moon system. It is associated with the right-handed sequence of three mutually orthogonal unit vectors $\{{\boldsymbol {{c}}_1}^{B}, \: {\boldsymbol {{c}}_2}^{B}, \: {\boldsymbol {{c}}_3}^{B}\}$, where ${\boldsymbol{c}}_1^B$ is aligned with the line joining Saturn and the moon at a reference epoch $t_0$ and points toward Saturn, ${\boldsymbol{c}}_2^B$ lies in the orbital plane of the system, and ${\boldsymbol{c}}_3^B$ is aligned with the orbital angular momentum vector.
The barycentric synodic frame ($\mathcal{R}_{SYN} ; {O_B; x, y, z}$) is a rotating frame that shares its origin with $\mathcal{N}_B$ and is identified by the right-hand triad $\{\boldsymbol {{i}}, \,  \boldsymbol {{j}}, \, \boldsymbol {{k}}\}$. Unit vector $\boldsymbol {{i}}$ points from the barycenter toward Saturn, and $\boldsymbol {{k}}$ is aligned with the angular momentum vector of the system. The synodic frame rotates with angular velocity $\boldsymbol{\omega}$ about $\boldsymbol{k} \equiv {\boldsymbol{c}}_3^{B}$, where $\omega$ denotes the mean motion of a Saturn–moon system. The rotation matrix from $\mathcal{N}_B$ to $\mathcal{R}_{SYN}$ is given by
\begin{equation}\label{rotatinframeeq}
    \underset{\mathcal{R}_{SYN} \leftarrow \mathcal{N}_B} {\rm \textbf{R}} = {\rm \textbf{R}}_3 \left( \omega \, \left( t - t_{0} \right) \right).
\end{equation}
At epoch $t_{0}$, the two frames coincide. {Fig. \ref{SYNframe} illustrates the two reference frames that share a common origin at the center of mass of the system.}

\indent The Radial–Transverse–Normal (RSW) reference frame ($\mathcal{R}_{RSW}; {O_{S/C}; \xi, \eta, \zeta}$) is a rotating coordinate system centered at the spacecraft and defined relative to the instantaneous orbital motion of the spacecraft. It is associated with the right-handed sequence of unit vectors $\{ \boldsymbol {{u}}_r, \; \boldsymbol {{u}}_\theta, \; \boldsymbol {{u}}_h \}$. In particular, $\boldsymbol {{u}}_r$ is aligned with the spacecraft position vector, $\boldsymbol {{u}}_h$ points toward the orbital angular momentum, and $\boldsymbol {{u}}_\theta$ completes the right-handed triad. {Fig.~\ref{rswframepic} shows the SCI frame $\mathcal{N}_S$ and the rotating RSW frame. }If $\Omega$, $i$, and $\theta_t$ denote the instantaneous longitude of the ascending node, inclination and argument of latitude, respectively, the rotation matrix from an inertial frame (i.e., SCI) to RSW is
\vspace{-4pt}
\begin{equation}\label{LVLHmat}
    \underset{\rm \mathcal{R}_{RSW} \leftarrow \mathcal{N}_S}{\rm \textbf{R}}
    =
     {\rm \textbf{R}}_3\left(\theta_t\right) {\rm \textbf{R}}_1\left(i\right) {\rm \textbf{R}}_3\left(\Omega\right).
\end{equation}

\subsection{Dynamical model for halo and science orbits}

\begin{figure}[h]
    \centering
    \includegraphics[width=0.5\textwidth]{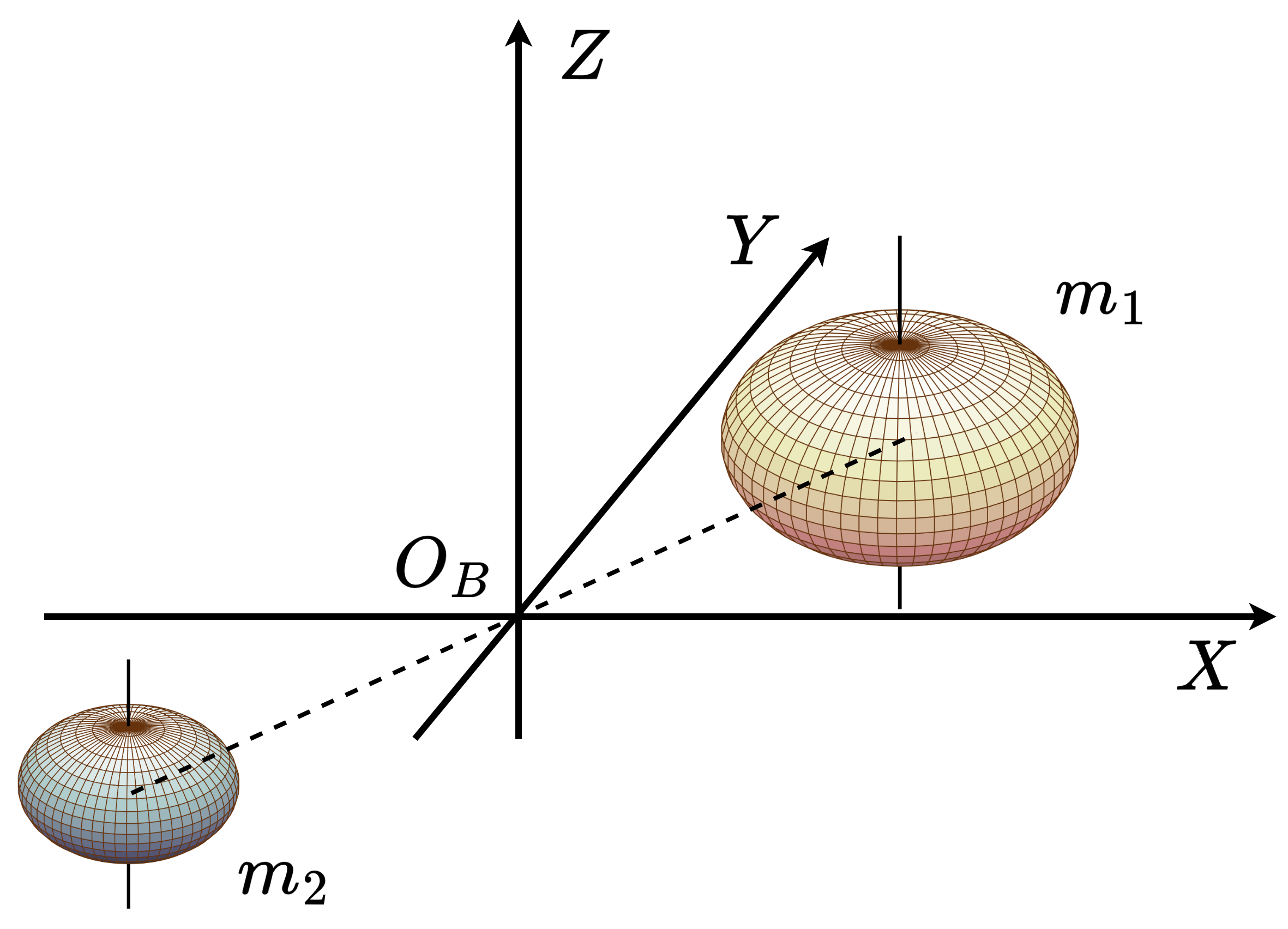}
    \caption{Scheme of a system consisting of two oblate, homogeneous spheroids with axes of maximum inertia parallel to each other. }
    \label{OblatePrimaries}
\end{figure}

Exploration orbits for lunar surface observations are modeled within the CR3BP framework. The system includes only Saturn and a single moon as massive bodies orbiting their common center of mass in circular orbits. This approximation is justified by the low eccentricities of the ILMs orbits (see Table \ref{tab:saturn_moons}). To enhance the fidelity of the model, the second zonal harmonic term ($J_2$) is introduced to account for the oblateness of the primaries. The resulting $J_2$–perturbed CR3BP considers two oblate spheroids of masses $m_1$ and $m_2$ separated by a constant distance $d$ (see Fig. \ref{OblatePrimaries}), with their axes of maximum inertia parallel to each other and to the orbital angular momentum of the system.\footnote{The alignment results from tidal locking, which gradually orients the axes of the bodies in the same direction, allowing them to be treated as parallel.} A third body of negligible mass (i.e., the spacecraft) is influenced only by the gravitational action of the primaries.

\begin{table*}[t]
\caption{Mass, equatorial radius $R$, mean orbital radius $d$, orbital period $T$, and orbit elements of the major moons of Saturn, together with the $J_2$ coefficients (raw data from \cite{spicenaif}). For the definition of the SOI radii, refer to Section~\ref{intermoonsec}.}
\centering
\small 
\begin{tabular}{@{}lllllllll@{}}
\toprule
\text{Moon} & \text{Mass} & \text{$R$} & \text{$d$} & \text{$T$} & \text{Inclination} & \text{Eccentricity} & \text{$J_2$} & \text{SOI radius} \\ 
& \text{[$10^{20}$ kg]} & \text{[km]} & \text{[$10^3$ km]} & \text{[day]} & \text{[degree]} & - & - & \text{[$10^3$ km]} \\\midrule
S    & $568.32 \cdot 10^{4}$  & 60268  & - & - & - & - & $1.6291 \cdot 10^{-2}$ & - \\ 
Mi      & 0.3799  & 207.4  & 185.52 & 0.9424 & 1.53 & 0.0202 & $3.1089 \cdot 10^{-2}$ & 1.123 \\ 
En & 1.0805   & 256.6  & 238.02 & 1.3702 & 0.01 & 0.0045 & $5.4352 \cdot 10^{-3}$ & 2.195 \\ 
Te   & 6.1742   & 540.4  & 294.66 & 1.8878 & 1.86 & 0.0000 & $9.4345 \cdot 10^{-3}$ & 5.463 \\ 
Di     & 10.954   & 563.8  & 377.40 & 2.7369 & 0.02 & 0.0022 & $1.4368 \cdot 10^{-3}$ & 8.791 \\ 
Rh      & 23.066   & 767.2  & 527.04 & 4.5175 & 0.35 & 0.0010 & $9.4961 \cdot 10^{-4}$ & 16.54 \\ 
Ti     & 1345.5 & 2575  & 1221.9 & 15.945 & 0.33 & 0.0292 & $3.3414 \cdot 10^{-5}$ & - \\ \bottomrule
\end{tabular}
\label{tab:saturn_moons}
\end{table*}

Let $X$, $Y$, and $Z$ denote the coordinates of the spacecraft in the $\mathcal{N}_B$ frame, centered at the center of mass of the system, with the equatorial planes of the primaries lying in the $XY$-plane (see Fig.~\ref{OblatePrimaries}). The gravitational potential of the system can be written as (see \cite{danby1992fundamentals})

\begin{equation}
U =
\frac{Gm_1}{r_1}
\left\{
1 - \frac{J_{2}^1}{2}
\left( \frac{R_1}{r_1} \right)^2
\!\left[3\!\left(\frac{Z}{r_1}\right)^2 - 1\right]
\right\}
+
\frac{Gm_2}{r_2}
\left\{
1 - \frac{J_{2}^2}{2}
\left( \frac{R_2}{r_2} \right)^2
\!\left[3\!\left(\frac{Z}{r_2}\right)^2 - 1\right]
\right\},
\label{eq:potential}
\end{equation}
where $G$ is the universal gravitational constant, $R_1$ and $R_2$ are the equatorial radii of the primaries, $r_1$ and $r_2$ denote the distances of the third body from $m_1$ and $m_2$, respectively, and $J_{2}^i$ ($i=1,2$) is the second zonal harmonic coefficient of each primary. 

Following the derivation in \cite{salazar2021science}, if $m_1$ and $m_2$ move in circular orbits about their common center of mass, the angular speed of the system is constant and equal to their orbital mean motion (see \cite{danby1992fundamentals})
\begin{equation}\label{meanmotionOblate}
   {
    n = \sqrt{\frac{G (m_1 + m_2)}{d^3} \left[ 1 + \frac{3}{2} \frac{(J_2^1 R_1^2 + J_2^2 R_2^2)}{d^2} \right]}
    },
\end{equation}

\noindent Changing to dimensionless variables, the total mass of the system $(m_1 + m_2)$ and the distance between the primaries $d$ are adopted as reference magnitudes, while time is scaled with respect to a reference value equal to
$\sqrt{\frac{d^3}{G \, (m_1 + m_2)}}$. This normalization yields a unit value for the gravitational constant. 

The normalized gravitational potential of the spacecraft in the $J_2$-perturbed system is

\begin{equation}\label{potentialJ2scFINAL}
    \tilde{U} =  {
    \cfrac{1 - \mu}{\tilde{r}_1}\left\{1 - \cfrac{A_1}{2 \tilde{r}_1^2}\left[3 \left(\cfrac{\tilde{Z}}{\tilde{r}_1}\right)^2 - 1\right] \right\} 
    } + {
    \cfrac{\mu}{\tilde{r}_2}\left\{1 - \cfrac{A_2}{2 \tilde{r}_2^2}\left[3 \left(\cfrac{\tilde{Z}}{\tilde{r}_2}\right)^2 - 1\right] \right\}.
    }
\end{equation}

\noindent where $\mu = \cfrac{m_2}{m_1 + m_2}$ and the variables marked with a tilde are $
    \tilde{r}_i = \cfrac{r_i}{d}, \quad \tilde{Z} = \cfrac{Z}{d}, \quad A_i = \cfrac{J_2^i \, R_i^2}{d}.
$

\noindent Using these normalized quantities, the mean motion in Eq. \eqref{meanmotionOblate} simplifies to
\begin{equation}\label{meanmotionJ2}
    \tilde{n} = \sqrt{1 + \frac{3 \, (A_1 + A_2)}{2}}.
\end{equation}

\begin{figure}[h] 
    \centering
\includegraphics[width=0.55\linewidth]{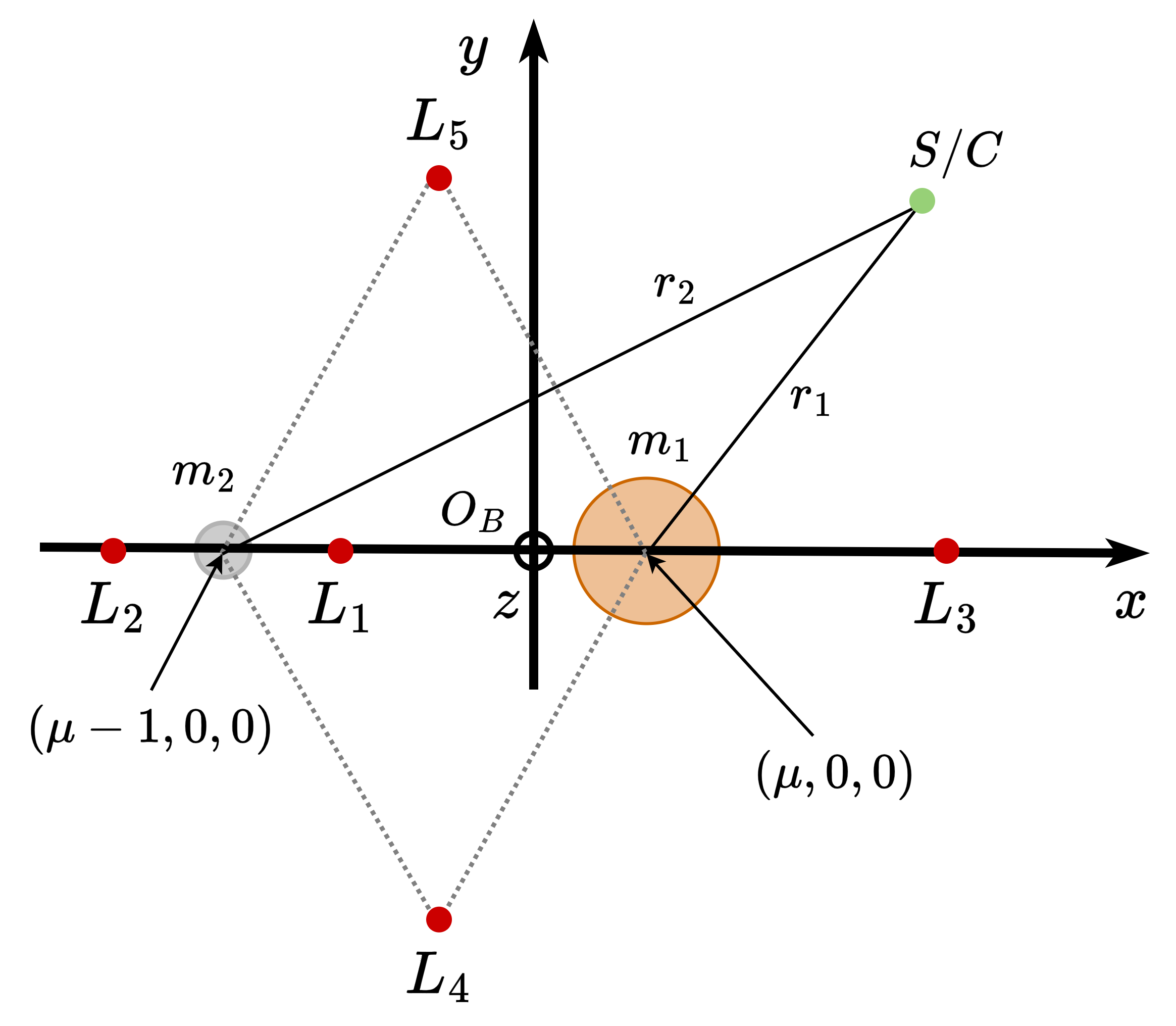}
    \caption{Configuration of the CR3BP showing the primaries ${m}_1$ and ${m}_2$, the spacecraft $S/C$, and the Lagrange equilibrium points in the synodic barycentric reference frame.}
    \label{CR3BP_SCHEME}
\end{figure}

\noindent Henceforward, normalized units are assumed, and the tilde will be omitted for simplicity. 
\noindent Consider now the barycentric synodic reference frame defined in Eq. \eqref{rotatinframeeq} rotating with angular velocity $n$ relative to the barycentric inertial frame, with the $x$–axis containing the larger primary at $(\mu,0,0)$  and the smaller one at $(\mu-1,0,0)$ (see Fig. \ref{CR3BP_SCHEME}). The equations of motion in this reference frame are

\begin{equation}
    \begin{cases}
       \ddot{x}-2n\dot{y} &= n^2x - (1-\mu)\cfrac{(x-\mu)}{r_1^3}C_1 - \mu\cfrac{(x+1-\mu)}{r_2^3}C_2,\\
\ddot{y}+2n\dot{x} &= y\!\left[n^2 - (1-\mu)\cfrac{C_1}{r_1^3} - \mu\cfrac{C_2}{r_2^3}\right],\\
\ddot{z} &= -z\!\left[(1-\mu)\cfrac{\tilde{C}_1}{r_1^3} + \mu\cfrac{\tilde{C}_2}{r_2^3}\right],  
    \end{cases}
\end{equation}

\noindent where
\begin{align}
C_i &= 1-\frac{3A_i}{2r_i^2}\!\left[5\!\left(\frac{z}{r_i}\right)^2-1\right], \quad i=1,2, \label{eq:C_i}\\
\tilde{C}_i &= C_i+\frac{3A_i}{r_i^2}, \quad i=1,2. \label{eq:Ctilde_i}
\end{align}

\noindent They can be rewritten in compact form 

\begin{equation}\label{eqMotJ2synNEW}
    \begin{cases}
        \ddot{x} &= \cfrac{\partial \Omega}{\partial x} + 2n\dot{y}, \\
    \ddot{y}  &= \cfrac{\partial \Omega}{\partial y} - 2n\dot{x},\\
    \ddot{z} &= \cfrac{\partial \Omega}{\partial z},
    \end{cases}
    \end{equation}
\noindent where $\Omega$ is the effective potential and is defined as follows:
\begin{equation}\label{OmegaPotential}
    \Omega = \,  U + \cfrac{n^2}{2}\left[ \left(1- \mu\right) \, r_1^2 + \mu \, r_2^2 - z^2\right] + \cfrac{\mu (1- \mu)}{2}.
\end{equation}

\noindent The dynamical system admits one integral of motion, the generalization of the Jacobi constant $C_J$ in the case of oblate primaries. The $J_2$-perturbed system admits five equilibrium points (Lagrangian points): three collinear points along the $x$-axis and two equilateral points (see \cite{salazar2021science} and Fig. \ref{CR3BP_SCHEME}). 


\subsection{Dynamical model for inter-moon transfers}
For the inter-moon phases of the tour, Saturn is modeled as the central body, and the spacecraft motion is described within a perturbed Saturn-spacecraft framework. The equations of motion are given by

\begin{equation}\label{2bp}
\begin{cases}
    \dot{\boldsymbol{r}} &= \boldsymbol{v} \\
    \dot{\boldsymbol{v}} &= -\cfrac{\mu_1}{r^3}\boldsymbol{r} + \cfrac{T}{m}\boldsymbol{u} + \boldsymbol{a}_p,\\
    \dot{m} &= -\cfrac{T}{g_0 I_\text{sp}},
\end{cases}
\end{equation}

\noindent where $\mu_1$ is the gravitational parameter of Saturn, $m$ is the mass of the spacecraft, $I_\text{sp}$ is the specific impulse of the thruster (assumed constant), $\boldsymbol{u}$ is the instantaneous direction of thrust, and $g_0 = 9.81 \, \mathrm{m/s^2}$. The state of the spacecraft, defined by position $\boldsymbol{r}$ and velocity $\boldsymbol{v}$, is expressed in the SCI frame. The term $\boldsymbol{a}_p$ represents the resultant of the non-Keplerian perturbative accelerations acting on the vehicle. To determine the perturbations to include in $\boldsymbol{a}_p$ and ensure that Eq. \eqref{2bp} accurately captures the spacecraft dynamics during inter-moon transfers, an evaluation of the perturbative effects is required. 

\section{Analysis and effects of orbit perturbations} \label{model}
Perturbation assessment is essential to streamline computations, in order to include only the significant contributions while excluding negligible ones, thus ensuring an efficient and accurate trajectory design process. {The purpose of this section is to quantify the effect of candidate perturbations on
inter-moon trajectories around Saturn following the approach described in \cite{flores2021method}, and to identify which of them must be retained in the dynamical model in Eq.  \eqref{2bp}.} Two main classes of perturbations are considered:
\begin{enumerate}[a)] 
 \item Third-body gravitational effects from the ILMs, Titan, Jupiter, and the Sun,
\item Oblateness effects associated with the $J_2$ terms of Saturn and its moons. 
\end{enumerate}


Let body $2$ represent the spacecraft, $1$ the primary body, and $3$ the perturbing body. When describing the relative motion of body 2 with respect to 1, the third-body gravitational perturbation due to 3 is
\begingroup
\thinmuskip=0mu
\medmuskip=0mu
\thickmuskip=0mu
\begin{equation}\label{a3B}
   {\boldsymbol{a}}_{3B} = \mu_3\,\left\{\frac{{\boldsymbol{r}}_{13} -{\boldsymbol{r}}_{12}}{\left[\left({\boldsymbol{r}}_{12} -{\boldsymbol{r}}_{13} \right)\cdot\left({\boldsymbol{r}}_{12} -{\boldsymbol{r}}_{13} \right)\right]^{\frac{3}{2}}} - \frac{{\boldsymbol{r}}_{13}}{r_{13}^3}\right\}.
\end{equation}
\endgroup

\noindent where $\mu_3$ is the gravitational parameter of the perturbing body, and ${\boldsymbol{r}}_{1j}$ ($j = 2, 3$) denotes the position vectors of body $j$ relative to body 1. {It is worth emphasizing that the positions of the moons and planets are obtained from the ephemeris and therefore already account for their mutual gravitational interactions. Moreover the mass of the spacecraft is negligible compared to the other bodies.}

For the perturbative acceleration due to the $J_2$ coefficient, let $\mu_j$ denote the gravitational parameter of the oblate body, $R_{j}$ its equatorial radius, $J_{2j}$ its second zonal harmonic coefficient, and $r$ the distance of the spacecraft from the body. The components of the perturbing acceleration due to the oblateness of the body are conveniently expressed in the RSW frame (see \cite{pontani2023advanced})

\begin{equation}\label{aJ2}
 {\boldsymbol{a}}_{J_2}^{RSW} = \frac{3 \mu_j}{r^4} R_{j}^2 J_{2_j}
\begin{bmatrix}
    \cfrac{3 \sin^2 (\omega + \theta) \sin^2 i - 1}{2} \\
    -\sin^2 i \sin (\omega + \theta) \cos (\omega + \theta) \\
    -\sin i \cos i \sin (\omega + \theta)
\end{bmatrix}^\text{T}
\begin{bmatrix}
    \boldsymbol {{u}}_r\\
    \boldsymbol {{u}}_\theta \\
   \boldsymbol {{u}}_h
\end{bmatrix}.
\end{equation}
\noindent Here, \(i\), \(\omega\), and \(\theta\) represent the inclination, the argument of periapsis, and the true anomaly of the spacecraft orbit with respect to the perturbing oblate body. Then, the resulting acceleration is transformed into the SCI frame, yielding ${\boldsymbol{a}}_{J_2}$. 

The total perturbative acceleration used in Eq. \eqref{2bp} becomes

\begin{equation}
    \boldsymbol{a}_P = \sum_{i=1}^{N} \boldsymbol{a}_{3B}^{(i)} +  \sum_{j=1}^{M} \boldsymbol{a}_{J_2}^{(j)},
\end{equation}
where $N$ is the total number of perturbing bodies considered for third-body effects, and $M$ is the total number of oblate bodies included for the $J_2$ perturbation.

A comparative analysis is conducted to quantify the contribution of each perturbation. For the moons, their instantaneous states are considered throughout the propagation, whereas Jupiter is assumed to be fixed at its minimum distance from Saturn, providing a conservative estimate of its influence. The initial state of the spacecraft lies on the equatorial plane of Saturn and is defined as

\begin{equation}\label{initialcondi}
\begin{split}
    \boldsymbol{r}_0 & = \left[x_0 \, \, \, 0 \, \, \, 0\right]^T \, \, \, \, \, \, \text{with} \, \, \, \, \, \, x_0 = \frac{r_{mA} + r_{mB}}{2}\\ 
    \boldsymbol{v}_0 & = \left[0 \, \, \, \dot{y}_0 \, \, \, 0\right]^T  \, \, \, \, \, \, \text{with} \, \, \, \, \, \, \dot{y}_0 = \sqrt{\frac{\mu_1}{\| \boldsymbol{r}_0\|}},
    \end{split}
\end{equation}
\noindent where $r_{mA}$ and $r_{mB}$ denote the mean Saturn-moon distances of the departing (A) and arrival (B) moons, respectively. The numerical integration of the equations of motion is carried out with relative and absolute tolerances of $\epsilon_{rel}=\epsilon_{abs}=10^{-13}$, determined from a convergence study to ensure the desired numerical accuracy \cite{flores2021method}. The propagation time span \(\Delta T_{\mathrm{syn}}\) is set to the maximum synodic period between the spacecraft unperturbed orbital period and the osculating periods of the two terminal moons (A and B). This ensures that the gravitational influence of both moons is fully captured during the propagation. 

\begingroup
The analysis is performed independently for each inter-moon leg (Rhea--Dione, Dione--Tethys, Tethys--Enceladus, and Enceladus--Mimas), and is carried out as follows:
\begin{enumerate}[1)]
   \item The initial state \((\boldsymbol{r}_0,\boldsymbol{v}_0)\) defined in Eq.~\eqref{initialcondi} is propagated using Eq.~\eqref{2bp} over \(\Delta T_{\mathrm{syn}}\), with the thrust set
    to zero (purely ballistic motion) and including all perturbations, yielding the reference final state
    \begin{equation}
        \boldsymbol{X}_f^{\,\text{(ref)}} = 
        \begin{bmatrix}
            \boldsymbol{r}_f^{\,\text{(ref)}} \\
            \boldsymbol{v}_f^{\,\text{(ref)}}
        \end{bmatrix}.
    \end{equation}

       \item Starting from the same initial state, the trajectory is propagated over \(\Delta T_{\mathrm{syn}}\) while omitting one perturbation at a time from the dynamical model. For each omitted perturbation \(i\), the final state is computed

    \begin{equation}
        \boldsymbol{X}_f^{(i)} =
        \begin{bmatrix}
            \boldsymbol{r}_f^{(i)} \\
            \boldsymbol{v}_f^{(i)}
        \end{bmatrix}.
    \end{equation}
    It is worth mentioning that, when the gravitational effect of a moon is neglected, the associated $J_2$ term is also excluded. 
    \item For each case, the final state \(\boldsymbol{X}_f^{(i)}\) is compared with  \(\boldsymbol{X}_f^{\,\text{(ref)}}\). The contribution of the omitted perturbation is quantified by the final position and velocity errors
    \begin{equation}
        e_r^{(i)} = \left\| \boldsymbol{r}_f^{(i)} - \boldsymbol{r}_f^{\,\text{(ref)}} \right\| ,
        \qquad
        e_v^{(i)} = \left\| \boldsymbol{v}_f^{(i)} - \boldsymbol{v}_f^{\,\text{(ref)}} \right\|.
    \end{equation}

    \item Following established literature \cite{salazar2021science,fantino2020design}, thresholds of 1 km and 1 m/s are set for $e_r^{(i)}$ and $e_v^{(i)}$, respectively. These values are chosen to meet the fidelity goals. If either error exceeds its corresponding threshold, the associated perturbation is classified as relevant and therefore included in the dynamical model for inter-moon transfers.
\end{enumerate}


\endgroup

\indent For the purpose of this analysis, the necessary parameters are retrieved from the JPL Navigation and Ancillary Information Facility (NAIF) archives \cite{spicenaif}. Specifically, the SPICE kernel \texttt{cpck15Dec2017.tpc}, which contains the most recent physical constants for the Saturn system as of December 15, 2017, is used. The SPICE toolkit, in conjunction with the aforementioned kernel, allows the extraction of physical parameters such as the gravitational constants, principal radii, and unnormalized \(J_2\) coefficients of the bodies of the Saturnian system (see Table \ref{tab:saturn_moons}). Additionally, the instantaneous state vectors of the perturbing bodies are obtained from the ephemeris provided by SPICE, expressed in the SCI frame. Although not critical to the perturbation analysis, the starting epoch for the propagation of the four scenarios was set to 2042-January-01, for consistency with previous work \cite{fantino2023end}.

Table \ref{perttab} summarizes the results of the analysis for the four scenarios. Contributions exceeding the defined thresholds are highlighted in red. The $J_2$ effect of Saturn is the dominant perturbation, with position errors of order $10^4$ km and velocity errors up to $10^3$ $[\nicefrac{\text{m}}{\text{s}}]$. The perturbing effects of the moons are also relevant and must be considered in trajectory design. The gravitational influence of the Sun appears to be significant only in the more external transfers, Rhea to Dione and Dione to Tethys, whereas the effect of Jupiter is negligible. Similarly, the $J_2$ perturbations of the moons produce errors well below the set thresholds. For completeness, the contribution of solar radiation pressure was also estimated. Assuming a fully reflecting surface, its effect in the most external segment (Rhea–Dione) is found to be an order of magnitude smaller than the already negligible gravitational perturbation from Jupiter. For this reason, solar radiation pressure is not included in the dynamical model.
Note that the gravitational perturbation of the rings of Saturn was also disregarded. This assumption is justified by Cassini mission data, which estimate the total ring mass as $(0.4 \pm 0.13) \, M_{\text{Mi}}$ (see \cite{iess2019measurement}), associated with negligible acceleration on Cassini and hence on the spacecraft considered in this study.\\

\begingroup 
\small  
\begin{table}[ht]
\centering
\caption{Results of the perturbation analysis for the four transfers. Values are presented in a logarithmic scale (e.g., $-1$ corresponds to an error of $10^{-1}$).}
{\renewcommand{\arraystretch}{1}
\begin{tabular}{l|cc|cc|cc|cc}
\toprule
\makecell[l]{{Perturbation} \\ {}} & 
\multicolumn{2}{c|}{Rh $\leftrightarrow$ Di} &
\multicolumn{2}{c|}{Di $\leftrightarrow$ Te} &
\multicolumn{2}{c|}{Te $\leftrightarrow$ En} &
\multicolumn{2}{c}{En $\leftrightarrow$ Mi} \\
\cmidrule(lr){2-3} \cmidrule(lr){4-5} \cmidrule(lr){6-7} \cmidrule(lr){8-9}
& \makecell{$e_r$ $[\text{km}]$} & \makecell{$e_v$ $[\nicefrac{\text{m}}{\text{s}}]$} & 
\makecell{$e_r$  $[\text{km}]$} & \makecell{$e_v$  $[\nicefrac{\text{m}}{\text{s}}]$} &
\makecell{$e_r$  $[\text{km}]$} & \makecell{$e_v$ $[\nicefrac{\text{m}}{\text{s}}]$} &
\makecell{$e_r$  $[\text{km}]$} & \makecell{$e_v$  $[\nicefrac{\text{m}}{\text{s}}]$} \\
\midrule
Mimas 3B        &  {\cellcolor{LightRed1}$\phantom{+}1$} & {${-1}$} & {\cellcolor{LightRed1}$\phantom{+}1$} & {${-1}$} & {\cellcolor{LightRed1}$\phantom{+}0$} & {${-1}$} & {\cellcolor{LightRed1}$\phantom{+}1$} & {\cellcolor{LightRed1}$\phantom{+}0$} \\
Enceladus 3B    & {\cellcolor{LightRed1}$\phantom{+}1$} & {${-1}$} & {\cellcolor{LightRed1}$\phantom{+}1$} & {${-1}$} & {\cellcolor{LightRed1}$\phantom{+}2$} & {\cellcolor{LightRed1}$\phantom{+}1$} & {\cellcolor{LightRed1}$\phantom{+}1$} & {\cellcolor{LightRed1}$\phantom{+}0$} \\
Tethys 3B       & {\cellcolor{LightRed1}$\phantom{+}2$} & {\cellcolor{LightRed1}$\phantom{+}0$}  & {\cellcolor{LightRed1}$\phantom{+}1$} & {\cellcolor{LightRed1}$\phantom{+}0$}    & {\cellcolor{LightRed1}$\phantom{+}2$} & {\cellcolor{LightRed1}$\phantom{+}1$} & {\cellcolor{LightRed1}$\phantom{+}2$} & {\cellcolor{LightRed1}$\phantom{+}0$} \\
Dione 3B        & {\cellcolor{LightRed1}$\phantom{+}3$} & {\cellcolor{LightRed1}$\phantom{+}0$}  & {\cellcolor{LightRed1}$\phantom{+}1$} & {\cellcolor{LightRed1}$\phantom{+}0$}    & {\cellcolor{LightRed1}$\phantom{+}1$} & {${-1}$} & {\cellcolor{LightRed1}$\phantom{+}2$} & {\cellcolor{LightRed1}$\phantom{+}0$} \\
Rhea 3B         & {\cellcolor{LightRed1}$\phantom{+}3$} & {\cellcolor{LightRed1}$\phantom{+}1$}  & {\cellcolor{LightRed1}$\phantom{+}2$} & {\cellcolor{LightRed1}$\phantom{+}1$}    & {\cellcolor{LightRed1}$\phantom{+}1$} & {${-1}$} & {\cellcolor{LightRed1}$\phantom{+}1$} & {\cellcolor{LightRed1}$\phantom{+}0$} \\
Titan 3B        & {\cellcolor{LightRed1}$\phantom{+}4$} & {\cellcolor{LightRed1}$\phantom{+}1$}  & {\cellcolor{LightRed1}$\phantom{+}2$} & {\cellcolor{LightRed1}$\phantom{+}1$}    & {\cellcolor{LightRed1}$\phantom{+}2$} & {\cellcolor{LightRed1}$\phantom{+}1$} & {\cellcolor{LightRed1}$\phantom{+}1$} & {\cellcolor{LightRed1}$\phantom{+}0$} \\
Jupiter 3B      & ${-1}$ & ${-3}$ & ${-1}$ & ${-3}$ & ${-3}$ & ${-4}$ & ${-2}$ & ${-3}$ \\
Sun 3B          & {\cellcolor{LightRed1}$\phantom{+}1$} & {\cellcolor{LightRed1}${0}$} & {\cellcolor{LightRed1}$\phantom{+}0$} & {\cellcolor{LightRed1}${0}$} & {${-1}$} & ${-2}$ & {${-1}$} & ${-2}$ \\
Saturn $J_2$    & {\cellcolor{LightRed1}$\phantom{+}4$} & {\cellcolor{LightRed1}$\phantom{+}2$} & {\cellcolor{LightRed1}$\phantom{+}4$} & {\cellcolor{LightRed1}$\phantom{+}2$}    & {\cellcolor{LightRed1}$\phantom{+}4$} & {\cellcolor{LightRed1}$\phantom{+}3$} & {\cellcolor{LightRed1}$\phantom{+}4$} & {\cellcolor{LightRed1}$\phantom{+}3$} \\
Mimas $J_2$     & ${-3}$ & ${-5}$ & ${-3}$ & ${-5}$ & ${-3}$ & ${-5}$ & ${-4}$ & ${-5}$ \\
Enceladus $J_2$ & ${-3}$ & ${-5}$ & ${-4}$ & ${-5}$ & ${-3}$ & ${-5}$ & ${-5}$ & ${-6}$ \\
Tethys $J_2$    & ${-3}$ & ${-5}$ & ${-3}$ & ${-4}$ & ${-3}$ & ${-5}$ & ${-4}$ & ${-5}$ \\
Dione $J_2$     & ${-3}$ & ${-5}$ & ${-3}$ & ${-4}$ & ${-4}$ & ${-5}$ & ${-4}$ & ${-5}$ \\
Rhea $J_2$      & ${-3}$ & ${-5}$ & ${-4}$ & ${-5}$ & ${-4}$ & ${-5}$ & ${-4}$ & ${-5}$ \\
Titan $J_2$     & ${-4}$ & ${-6}$ & ${-3}$ & ${-5}$ & ${-3}$ & ${-5}$ & ${-4}$ & ${-5}$ \\
\bottomrule
\end{tabular}
}
\label{perttab}
\end{table}
\endgroup

\section{Tour design strategy}\label{tourdesign}
This section outlines the methodology adopted to design the lunar tour by leveraging the natural dynamics of the system in combination with LT propulsion.
Halo orbits and their associated dynamical structures are presented and used to build exploration orbits around individual moons as well as moon-to-moon transfer trajectories.
The overall design strategy is divided into two main components: science orbits definition and inter-moon transfers planning.  

\subsection{Science orbits design}\label{sciencedesign}
Dynamical substitutes of the unperturbed halo orbits exist in the $J_2$-perturbed model and serve as staging and science orbits near the Lagrangian points $L_1$ and $L_2$ of the Saturn-moon systems (whose $x$-coordinates are listed in Table~\ref{LagPoints}) \cite{danby1992fundamentals}. These orbits are organized into families, parameterized by their energy levels represented by the Jacobi constant. 

{\small
\begin{table}[h]
\caption{Mass ratios $\mu$ of the Saturn–moon systems and corresponding $x$-coordinates of the $L_1$ and $L_2$ equilibrium points in the barycentric synodic frame of the $J_2$–perturbed CR3BP (in the respective distance units).}
\vspace{-0.2cm}
\centering
\begin{tabular}{lccc}
\hline
{System} & $\mu$ & $x_{L_1}$ & $x_{L_2}$ \\
\hline
SMi & $6.67 \times10^{-8}$ & $-0.9971878$ & $-1.0028172$\\
SEn & $1.89 \times10^{-7}$ & $-0.9960222$ & $-1.0039880$\\
STe & $1.09\times10^{-6}$ & $-0.9928877$ & $-1.0071440$ \\
SDi & $1.93\times10^{-6}$ & $-0.9913849$ & $-1.0086609$ \\
SRh & $4.06\times10^{-6}$ & $-0.9889734$ & $-1.0111001$ \\
\hline
\end{tabular}
\label{LagPoints}
\end{table}}

\begin{figure*}[h]
    \centering
    \begin{subfigure}[b]{0.48\textwidth}
        \centering
        \includegraphics[width=\textwidth]{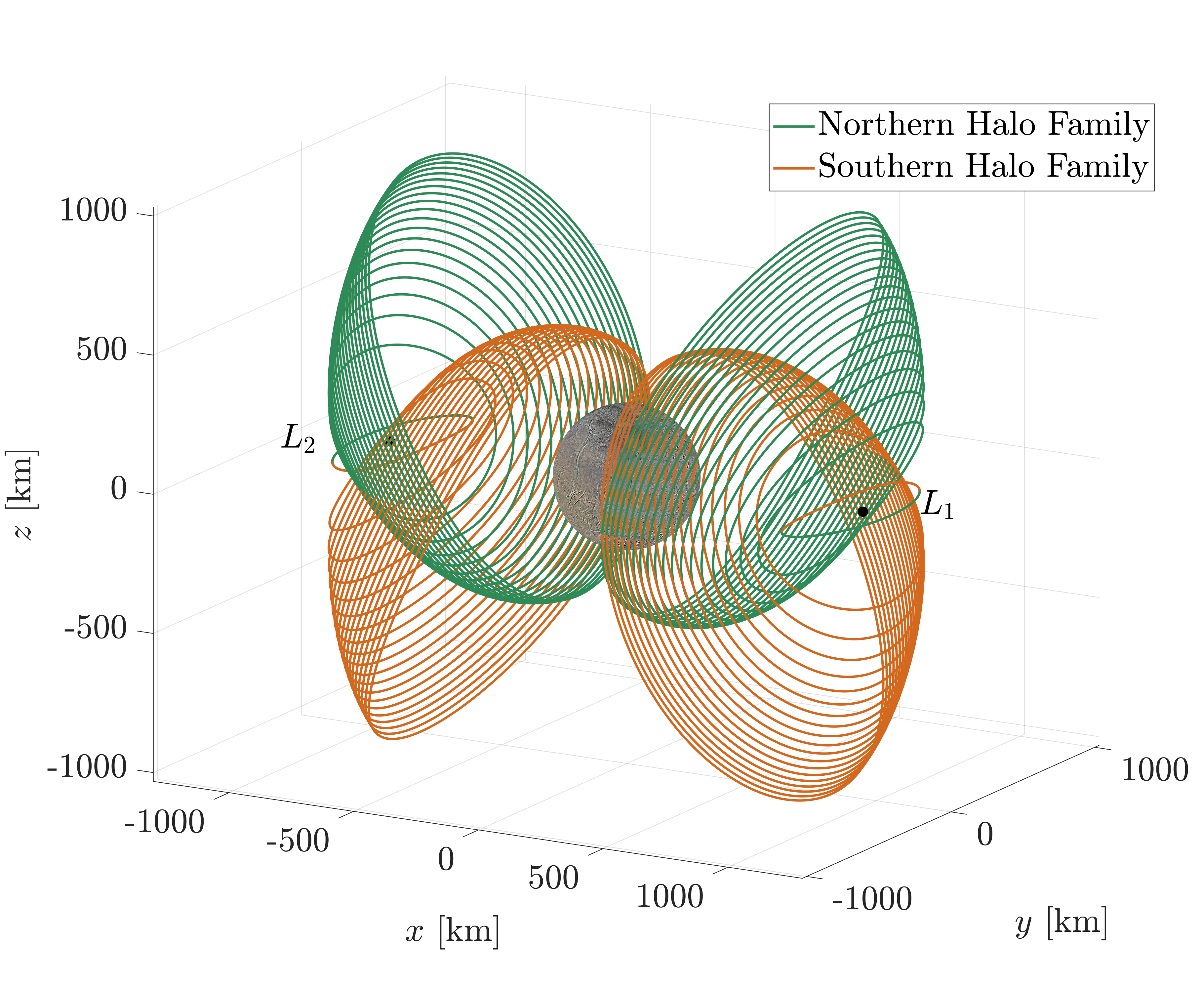}
        \caption{Three-dimensional side view.}
    \end{subfigure}
    \hfill
    \begin{subfigure}[b]{0.48\textwidth}
        \centering\includegraphics[width=\textwidth]{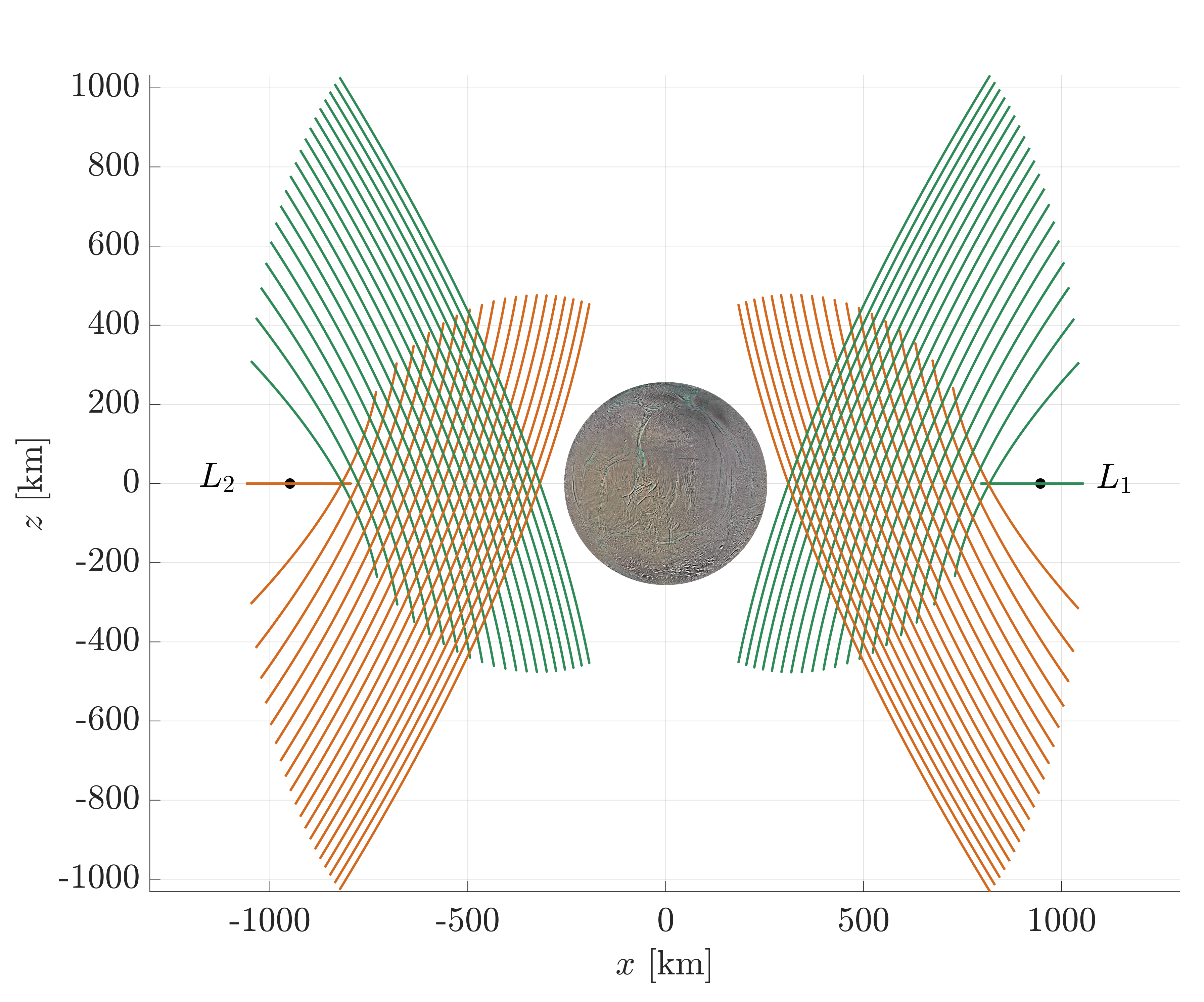}
        \caption{Projection onto the \(xz\)-plane.}
    \end{subfigure}
    \caption{Families of halo orbits around $L_1$ and $L_2$ in the $J_2$-perturbed Saturn–Enceladus CR3BP (Enceladus-centered synodic reference frame). The families are characterized by $C_J \in [3.00252, 3.00257]$, with a $C_J$ spacing of $2.0\times10^{-6}$ between consecutive orbits.}
    \label{fig:HaloEnc}
\end{figure*}

Halo orbits are generated using a differential-correction scheme applied to the initial conditions of the unperturbed model, enforcing $xz$-plane symmetry and periodicity. The numerical algorithm described in~\cite{connor1984three} is employed for this purpose. Each libration point admits two families: Northern halos, extending mostly above the \(xy\)-plane (positive $z$), and Southern halos, which lie predominantly below it (negative $z$). The families are computed using continuation methods, with a uniform discretization in $C_J$ and a fine step along the $x$-coordinate of the initial state. Figure \ref{fig:HaloEnc} shows an example of the Northern and Southern halo families around $L_1$ and $L_2$ for the Saturn–Enceladus system.

\begin{figure}
    \centering
    \includegraphics[width=0.55\linewidth]{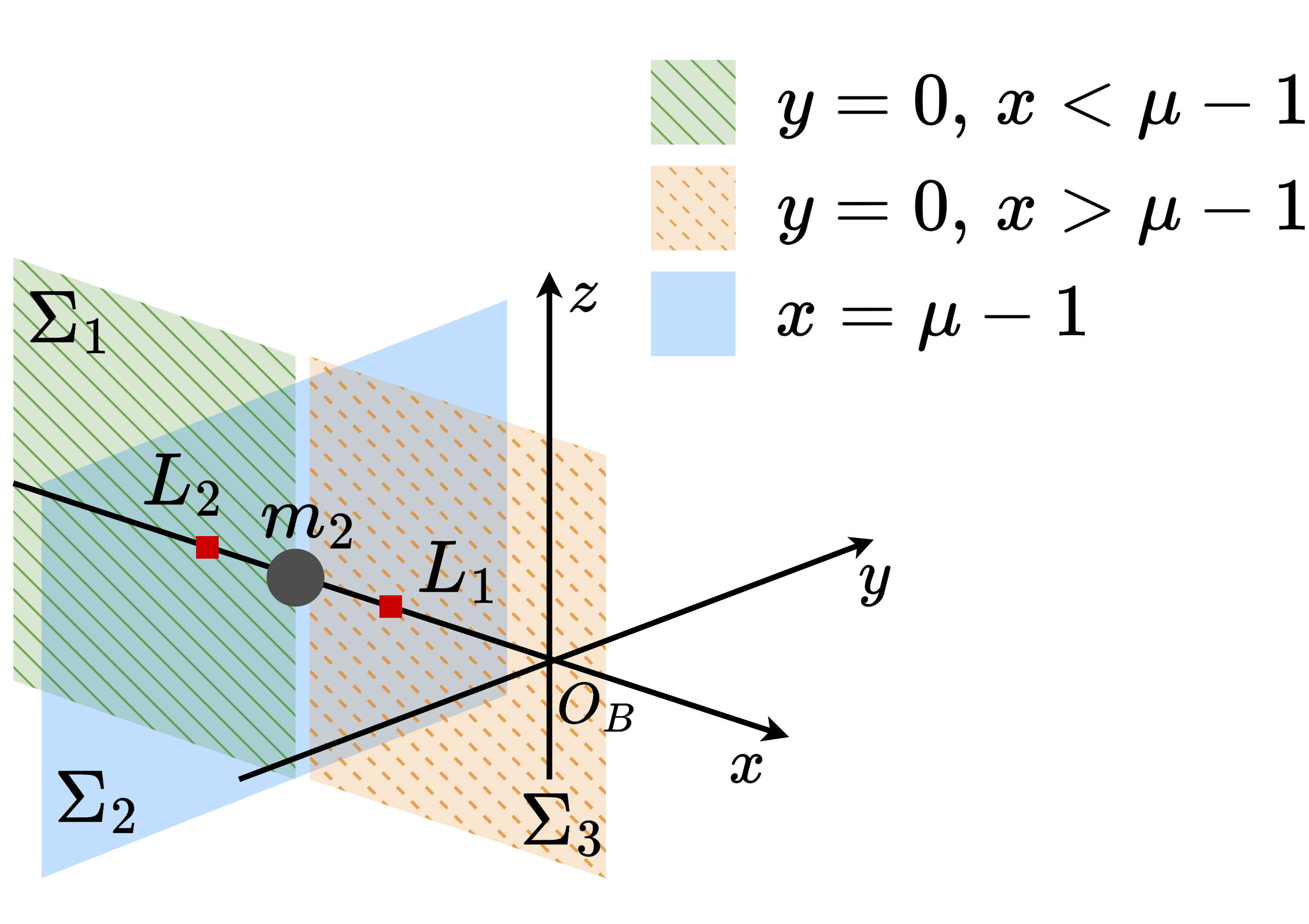}
     \caption{Representation of the Poincaré sections used in the analysis.}
    \label{fig:PsectionPic}
\end{figure}

\indent Dynamical structures of the CR3BP persist in the $J_2$-perturbed framework, including stable and unstable HIMs associated with the halo orbits. HIMs are computed and numerically propagated using established methods (see \cite{parker2012practical}). A small perturbation of $10^{-6}$ in normalized units is applied to points along the periodic orbit, in the directions of the stable and unstable eigenvectors of the monodromy matrix, to define the initial states of the associated manifolds.
\noindent Each halo orbit is discretized into a finite number of points, each point serving as the origin of an invariant manifold trajectory. \\
\indent The HIMs enable the design of low-energy transfers and fuel-efficient science orbits that provide global coverage of the moons. This part of the work focuses on designing homoclinic and heteroclinic connections. For a periodic orbit $\Gamma$, a homoclinic connection is a trajectory that asymptotically approaches $\Gamma$ both forward and backward in time, and is obtained by the intersection of its unstable and stable manifolds. A heteroclinic connection between two periodic orbits, $\Gamma$ and $\Lambda$, is a trajectory that tends forward in time to $\Lambda$ and backward in time to $\Gamma$, or vice versa. It is identified by the intersection of the unstable manifold of $\Gamma$ and the stable one of $\Lambda$.\\
\indent Manifold connections are obtained by propagating the HIMs until their first or second crossing (depending on the connection type) of a Poincaré section, a predefined plane \(\Sigma\) in phase space (see \cite{fantino2020design}). For homoclinic and heteroclinic connections between halo families around the same libration point, two Poincaré sections are employed: 
$\Sigma_1$, defined by $y = 0$ and $x < \mu - 1$, is used for connections departing from and returning to orbits around $L_1$;  
$\Sigma_3$, defined by $y = 0$ and $x > \mu - 1$, is used for connections departing from and returning to orbits around $L_2$.
For heteroclinic connections between different libration points, the plane $\Sigma_2$ is $x = \mu - 1$, orthogonal to the $x$–axis and passing through the center of the moon. Figure~\ref{fig:PsectionPic} illustrates the geometry of these three Poincaré sections. A minimum safety distance of 20 km from the surface of the moons is maintained.  A zero-cost connection is achieved when the states of the stable and unstable manifolds coincide in magnitude and direction at \(\Sigma\). Thus, identifying the best connection involves finding the pair of unstable and stable trajectories with the minimum position and velocity differences at the section. Tolerances of \( 1 \, \text{km} \) for positional error and \( 1 \, \text{m/s} \) for velocity error are adopted, consistent with the thresholds used in Section~\ref{model} and with values employed in literature \cite{salazar2021science,fantino2020design}. A notable advantage of this method is that the heteroclinic and homoclinic connections can be swept for an arbitrary duration with minimal to negligible propellant consumption. Their looping nature supports comprehensive observational campaigns of the lunar surfaces, including polar regions, thus improving the scientific return of the mission.

\indent To explore all potential combinations of departure and arrival orbits, the following heteroclinic connection types are considered:
\begin{description}
    \item[(A)] Northern \(L_1\) to Northern \(L_2\)
    \item[(B)] Northern \(L_1\) to Southern \(L_2\)
    \item[(C)] Northern \(L_1\) to Southern \(L_1\)
    \item[(D)] Southern \(L_2\) to Northern \(L_2\)
\end{description}

\noindent Additionally, the following homoclinic connections are investigated:
\begin{description}
    \item[(E)] Northern \(L_1\) to Northern \(L_1\)
    \item[(F)] Southern \(L_2\) to Southern \(L_2\)
\end{description}
\noindent For completeness, the inverse transfers have been computed, but they are omitted here for brevity based on symmetry considerations. For heteroclinic connections between halo orbits around the same equilibrium point, the search was restricted to pairs of halos with the same \(C_J\) values. 
For connections involving halos around different equilibrium points, each periodic orbit is associated with the one in the opposite family with the closest value of $C_J$. Then connections are sought between these minimum-$\Delta C_J$ pairs. Due to the geometric characteristics of the HIMs, Type A trajectories are identified at the second intersection of the HIMs with $\Sigma_2$, while Type B connections are determined at the first intersection with the same plane. Types~C, D, E, and F result from the first intersection with $\Sigma_1$ or $\Sigma_3$, as defined above.

\subsection{Inter-moon transfers} \label{intermoonsec}
Invariant manifolds associated with halo orbits are a powerful tool for designing low-cost transfers between moons.  However, as detailed in \cite{fantino2019connecting},
the HIMs do not intersect in configuration space due to the near-circular nature of the Saturn-centered orbits, necessitating orbital maneuvers for inter-moon transfers. To address this, a LT control law that minimizes propellant consumption through an electric propulsion system is implemented. The spatial outward branches of the HIMs are numerically propagated until their intersection with the SOI of the terminal moons. The SOI radius is defined as the Laplace SOI for the moon, scaled by an ad hoc factor \( f \): 
\begin{equation}
    R_{SOI} = f \, d \, \left({\frac{m_2}{m_1}} \right)^{\frac{2}{5}}.
\end{equation}
\noindent where $d$ is the distance between the two primaries. The factor \( f \) is chosen to be large enough to ensure that the SOI contains all halo orbits of both families and intersects the global manifold transversely. In this work, \( f \) is set to 4.5 for all five Saturn–moon systems, and the resulting SOI radii are listed in Table \ref{tab:saturn_moons}. Each transfer originates from a halo orbit around the \( L_1 \) point of the outer moon and terminates at one around \( L_2 \) of the inner moon. Within the SOIs, unstable manifolds are leveraged for the outer (departure) moon and stable ones for the inner (arrival) moon. The state vectors of the HIMs at the SOI are transformed in the SCI frame, and the elements of the corresponding osculating 2BP orbits with focus at Saturn are computed. \\
\begin{figure}[h]
    \centering
    \includegraphics[width=0.45\linewidth]{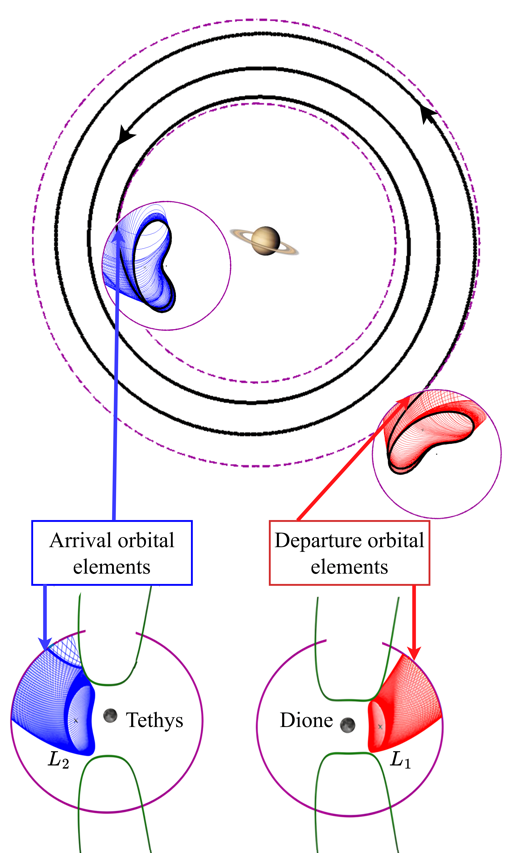}
    \caption{Illustration of the concept for consecutive moon-to-moon transfers.}
    \label{moontourex}
\end{figure}
\indent The propelled arc connecting the states at the two SOIs is modeled as a continuous LT transfer, illustrated by the black spiraling curve in Fig. \ref{moontourex}. 
\begingroup
The propulsion system characteristics (the maximum thrust $T_{\max}$, the
specific impulse $I_{sp}$, and the input power $P_{\text{in}}$) are summarized in
Table~\ref{tab:propulsion}. In this study, a constant thrust magnitude equal to
$T_{\max}$ is assumed, consistent with the performance of the
PPS\textsuperscript{\textregistered}X00 Hall-effect thruster, with $P_{\text{in}}$ provided by RTGs. The selected parameters are
compatible with existing electric propulsion systems and reflect a trade-off
between power consumption, thrust level, and mission duration.

\begin{table}[t]
    \centering
    \small
    \caption{Assumed thruster performance parameters.}
    \label{tab:propulsion}
    {
    \begin{tabular}{c c c}
        \hline
        $T_{\max}$ & $I_{sp}$ & $P_{\text{in}}$ \\
        \hline
        36 mN & 1600 s & 640 W \\
        \hline
    \end{tabular}}
\end{table}
\endgroup

The thrust direction is governed by a guidance law that locally maximizes the rate of reduction of an error function $Q$, defined in terms of osculating orbital elements $\mathbf{o} =\left[
a \, \, e \, \, i \, \, \omega \right]^T$, (semi-major axis, eccentricity, inclination, and argument of periapsis) as follows:
\begin{equation}\label{errfun}
\begin{split}
Q = & w_1 \left(a - a^*\right)^2 + w_2 \left(e \cdot a - e^* \cdot a^*\right)^2 + \\ 
+ & w_3  \left(i - i^*\right)^2 + 
w_4 \left\{e \cdot \cos^{-1}\left[\cos(\omega - \omega^*)\right]\right\}^2,
\end{split}
\end{equation}

\noindent where the superscript $*$ denotes the corresponding orbital elements of the target orbit. The coefficients $w_i > 0 $ $(i=1,\dots,4)$ are scalar weights that must be tuned for convergence and performance. {In particular, they are optimized, with the objective of minimizing the transfer time, while driving the orbit elements toward the desired values, with prescribed tolerances. It has been shown that, for an appropriately tuned set of weights, the guidance-law
formulation yields LT orbit transfers that approach minimum-time ones \cite{petropoulos2005optimisation}. It is worth emphasizing that the time-optimal transfers are also the most efficient option in terms of propellant consumption, because the thrust magnitude is constant and equal to $T_{max}$.} In Eq. \eqref{errfun} the eccentricity is multiplied by $a$ to maintain dimensional homogeneity. The term involving $\omega$ measures the shortest angular distance between two positions on a circle, restricted to the interval $[0,\,\pi]$. The longitude of the ascending node, $\Omega$, is excluded from the targeted element set, as its adjustment corresponds to a phase correction in true anomaly $\theta$, provided that the orbit of the arrival moon is circular and equatorial. This assumption holds for the Saturnian system, whose moons exhibit very low eccentricities and inclinations. Moreover, the time scale of the LT transfers is typically at least one order of magnitude larger than the orbital periods of the moons, and this allows disregarding phasing requirements. Phase corrections can be introduced post hoc by fine-tuning the thrust level. At this stage of the analysis, targeting the four orbital elements \((a, \, e, \, i, \, \omega)\) is sufficient to prove mission feasibility. The error function (Eq. \eqref{errfun})  can be regarded as a specific application of the Proximity Quotient guidance law (Q-Law) \cite{petropoulos2004low}.
\begingroup
The adopted strategy consists in maximizing the instantaneous rate of reduction of the error function $Q$. The thrust acceleration vector is defined as
\begin{equation}
\label{thrustacc}
\mathbf{f} =
\frac{T_{\max}}{m} \,\mathbf{u},
\qquad \|\mathbf{u}\| = 1,
\end{equation}
where $\mathbf{u}$ is the instantaneous thrust direction in the RSW frame.
Since the magnitude of $\mathbf{f}$ is small compared to the gravitational
pull of Saturn, its effect can be linearized (i.e., the rate of change of
$Q$ is a homogeneous function of the thrust components). Under the assumption that the propulsion system is always throttled at maximum thrust,
the only control variable is $\mathbf{u}$. Then the time derivative of the
error function $Q$ follows from the chain rule:
\begin{equation}\label{qdot}
\dot Q =
\nabla Q^{\mathsf{T}}\frac{\mathrm{d}\mathbf{o}}{\mathrm{d}t}
=
\nabla Q^{\mathsf{T}} \mathbf{R} \mathbf{f}.
\end{equation}
where $\nabla Q$ is the gradient of $Q$ with respect to the four targeted
orbital elements, and $\mathrm{d}\mathbf{o}/\mathrm{d}t = \mathbf{R}\mathbf{f}$
are the respective Gauss' planetary equations in matrix form \cite{pontani2023advanced}. Here $\mathbf{R}$ is a $4\times 3$ matrix whose entries are the partial derivatives of the orbital element rates with respect to the thrust acceleration
components
\begin{equation}
\mathbf{R}
= \frac{\partial \dot{\mathbf{o}}}{\partial \mathbf{f}}
= \left[\,
\frac{\partial \dot{o}_k}{\partial f_j}
\,\right]_{k=1,\dots,4}^{\;\;j=r,\theta,h}.
\end{equation}
Substituting Eq.~\eqref{thrustacc} into Eq.~\eqref{qdot} yields
\begin{equation}
\dot Q = \frac{T_{\max}}{m}\bigl(\mathbf{R}^{\mathsf{T}}\nabla Q\bigr)^{\mathsf{T}} \mathbf{u}.
\end{equation}
Because the thrust magnitude is fixed, minimizing $Q$ in the shortest time is equivalent to making $\dot Q$ as negative as possible at each instant. This is achieved by aligning the thrust direction opposite to $\mathbf{R}^{\mathsf{T}}\nabla Q$. Therefore, the optimal control law in the RSW frame is
\begin{equation}
\mathbf{f}^*
=
-\frac{T_{\max}}{m}\,
\frac{\mathbf{R}^{\mathsf{T}} \nabla Q}
{\bigl\|\mathbf{R}^{\mathsf{T}} \nabla Q\bigr\|}.
\end{equation}
Then the thrust acceleration vector is transformed into the SCI frame and included in Eq.~\eqref{2bp} for numerical propagation.
\endgroup

The use of LT propulsion offers improved fuel efficiency, but it results in longer flight times compared to direct impulsive transfers using chemical propulsion. {Based on the analysis presented in Section \ref{model}, the propagation of each segment includes all the perturbations identified as relevant. The instantaneous positions and velocities of the perturbing bodies are obtained from SPICE ephemeris, ensuring a high level of dynamical fidelity.}


\section{Results} \label{resultssec}
This section presents the results of the work, organized in two parts. The first focuses on the science orbits obtained within the perturbed CR3BP, describing their characteristics and observational properties. The second part discusses the results of the inter-moon transfer phases of the tour, reporting the trajectories and their associated performance metrics.

\subsection{Science orbits}\label{resultsSO}
For each moon, families of halo orbits are computed over a broad range of Jacobi constants, $C_J$, using an energy discretization such that each of the four families contains 25 members. 
The corresponding ranges of $C_J$ and orbital periods for the five $J_2$–perturbed Saturn–moon CR3BPs are reported in Table~\ref{tab:halo_parameters}.

\begin{table}[H]
\caption{Jacobi constant and period ranges for halo orbits in the Saturn-moon systems.}
\centering
{\renewcommand{\arraystretch}{1.1}
\begin{tabular}{lll}
\hline
{System} & {$C_J$ range} & {Period [hour]} \\
\hline
SMi & $[3.00423,\, 3.00425]$ & $[10.01,\, 11.08]$ \\
SEn & $[3.00252,\, 3.00257]$ & $[14.49,\, 16.07]$ \\
STe & $[3.00137,\, 3.00157]$ & $[17.14,\, 22.13]$ \\
SDi & $[3.00054,\, 3.00086]$ & $[21.72,\, 32.01]$ \\
SRh & $[2.99975,\, 3.00022]$ & $[31.73,\, 52.75]$ \\
\hline
Orbit Index & \hspace{0.1cm}1 \hspace{1cm} 25 & \hspace{0.1cm} 25 \hspace{0.6cm} 1 \\
\hline
\end{tabular}}
\label{tab:halo_parameters}
\end{table}

\begin{table}[ht]
\centering
\small
\caption{Characteristics of the heteroclinic connections between halo orbits in the five Saturn-moon systems.}
\label{tabresults_connections1}
\setlength{\tabcolsep}{9pt}
\renewcommand{\arraystretch}{0.8}
\begin{tabular}{ll|c|c|c|c|c}
\toprule
 & & {SMi} & {SEn} & {STe} & {SDi} & {SRh} \\
\midrule
\multirow{4}{*}{\rotatebox{90}{{Type A}}}
& $C_J$             & 3.004255 & 3.002566 & 3.001562& 3.000832 & 3.000166\\
& $\Delta T$ [hour] & 32.83 & 49.64 & 69.41 & 102.6 & 178.6\\
& $\Delta r$  [km]  & 0.51 & 0.37 & 0.13 & 0.35 & 0.51\\
& $\Delta v$ [m/s]  & 5.22 & 1.95 & 12.1 & 11.3& 9.84\\
\midrule
\multirow{4}{*}{\rotatebox{90}{{Type B}}}
& $C_J$             & - & 3.002566 & 3.001553 & 3.000821 & 3.000173 \\
& $\Delta T$ [hour] & - & 39.32& 57.72 & 85.34 & 144.7 \\
& $\Delta r$  [km]  & - & 0.27 & 0.62 & 0.53 & 0.73 \\
& $\Delta v$ [m/s]  & - & 0.61& 5.42 & 6.42 & 12.01\\
\midrule
\multirow{4}{*}{\rotatebox{90}{{Type C}}}
& $C_J$             & 3.004242 & 3.002540 & 3.001467 & 3.000694 & 3.000225\\
& $\Delta T$ [hour] & 38.56 & 58.36 & 83.96 & 123.6& 146.9\\
& $\Delta r$  [km]  & 0.49 & 0.93 & 0.86 & 0.92 & 0.81\\
& $\Delta v$ [m/s]  & 0.47 & 6.62 & 1.53 & 2.79 & 14.5\\
\midrule
\multirow{4}{*}{\rotatebox{90}{{Type D}}}
& $C_J$             &  3.004243 & 3.002539 & 3.001485 & 3.000694 & 2.999972 \\
& $\Delta T$ [hour] & 36.47 & 56.92 & 81.24 & 125.8 & 147.1\\
& $\Delta r$  [km]  & 0.27 & 0.88 & 0.58 & 0.24& 2.55\\
& $\Delta v$ [m/s]  & 0.42 & 0.13 & 3.15 & 8.95& 0.67\\
\bottomrule
\end{tabular}
\medskip
\footnotesize
\end{table}

\begin{table}[h]
\centering
\small
\caption{Characteristics of the homoclinic connections between halo orbits in the five Saturn-moon systems.}
\label{tabresults_connections2}
\setlength{\tabcolsep}{8pt}
\renewcommand{\arraystretch}{0.8}
\begin{tabular}{ll|c|c|c|c|c}
\toprule
 & & {SMi} & {SEn} & {STe} & {SDi} & {SRh} \\
\midrule
\multirow{4}{*}{\rotatebox{90}{{Type E}}}
& $C_J$             & - & 3.002573 & 3.001575 & 3.000853 & 3.000223\\
& $\Delta T$ [hour] & - & 48.85 & 71.87 & 105.9 & 146.9\\
& $\Delta r$  [km]  & - & 0.14 & 0.69 & 0.95 & 0.81\\
& $\Delta v$ [m/s]  & - & 0.23 & 1.27 & 3.10 & 13.6\\
\midrule
\multirow{4}{*}{\rotatebox{90}{{Type F}}}
& $C_J$             & - & 3.002573 & 3.001577 & 3.000856 & 3.000223\\
& $\Delta T$ [hour] & - & 49.31 & 72.83 & 107.6 & 147.5\\
& $\Delta r$  [km]  & - & 0.11 & 0.83 & 0.27 & 0.38\\
& $\Delta v$ [m/s]  & - & 0.13 & 1.32 & 2.42 & 6.17\\
\bottomrule
\end{tabular}
\medskip
\footnotesize
\end{table}

\begin{figure*}[htpb]
    \centering
    \begin{subfigure}{0.65\textwidth}
        \centering
        \includegraphics[width=\linewidth]{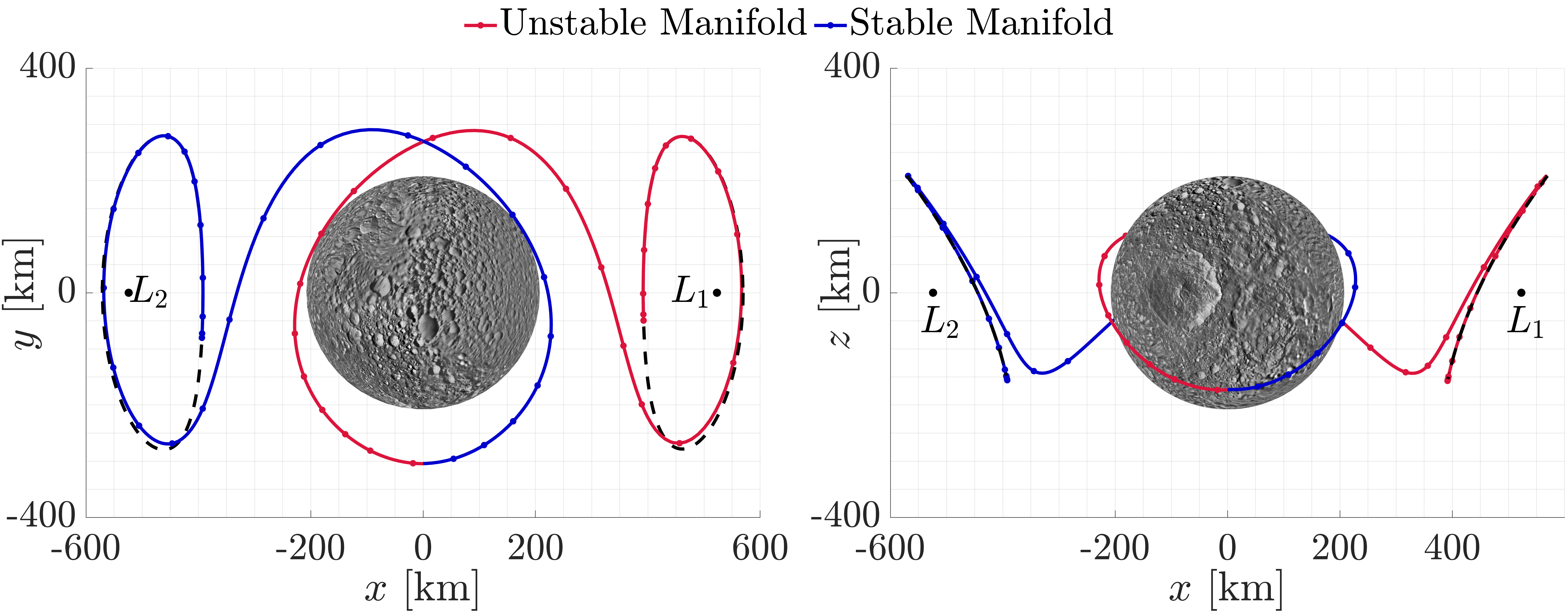}
        \caption{Type A - Saturn-Mimas system.}
    \end{subfigure}
    \vspace{0.4cm}
    \begin{subfigure}{0.65\textwidth}
        \centering
        \includegraphics[width=\linewidth]{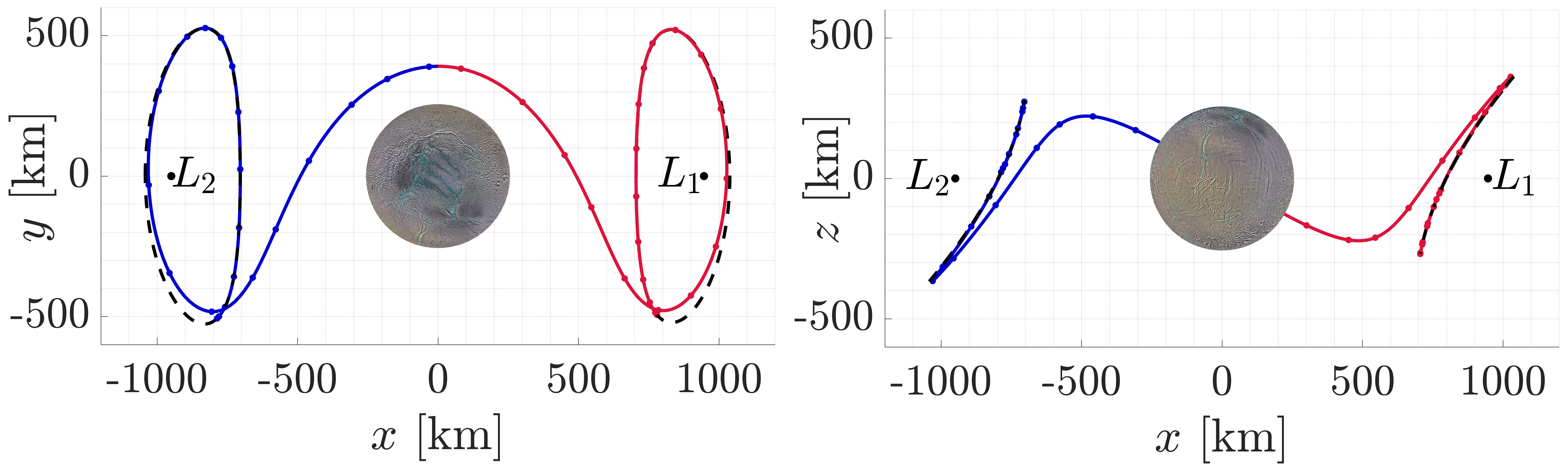}
                \caption{Type B - Saturn-Enceladus system.}
    \end{subfigure}
    \vspace{0.4cm}
    
    \begin{subfigure}{0.6\textwidth}
        \centering
        \includegraphics[width=\linewidth]{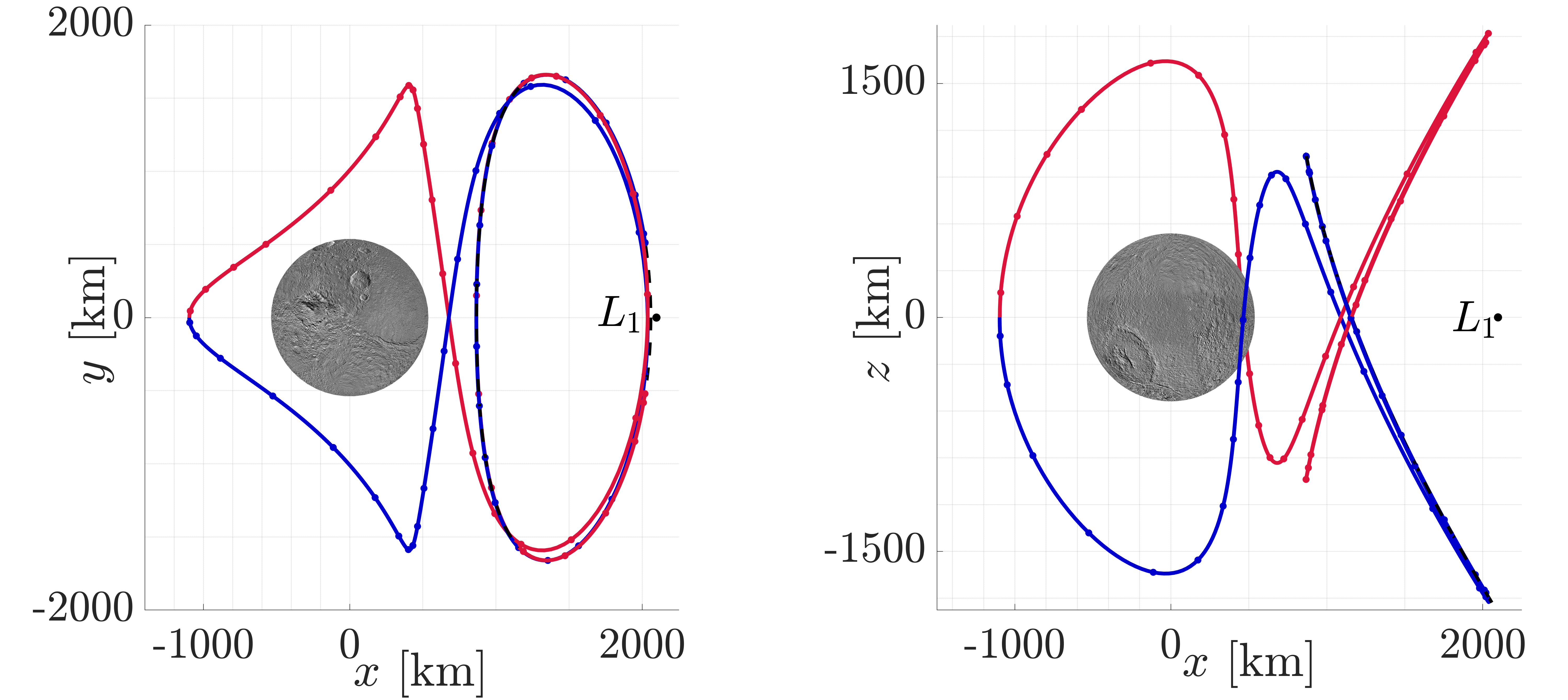}
                \caption{Type C - Saturn-Tethys system.}
    \end{subfigure}
    \vspace{0.4cm}
    \begin{subfigure}{0.6\textwidth}
        \centering
        \includegraphics[width=\linewidth]{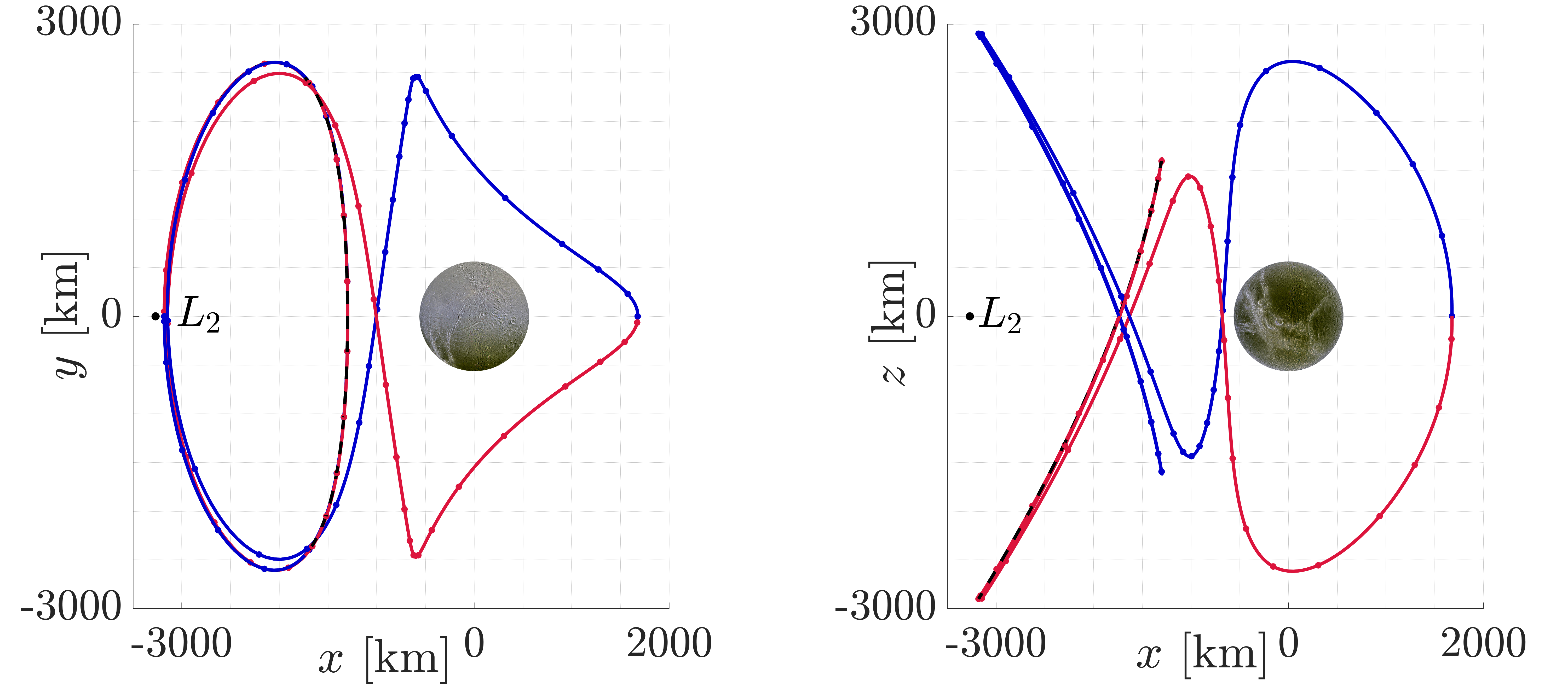}
                \caption{Type D - Saturn-Dione system.}
    \end{subfigure}
        \vspace{-0.2cm}
    \caption{Planar projections of representative heteroclinic connections (Types A to D) in the $J_2$-perturbed CR3BP (moon-centered synodic reference frame).}
    \label{fig:connections1}
\end{figure*}

\begin{figure*}[h]
    \centering
    \includegraphics[width=0.7\linewidth]{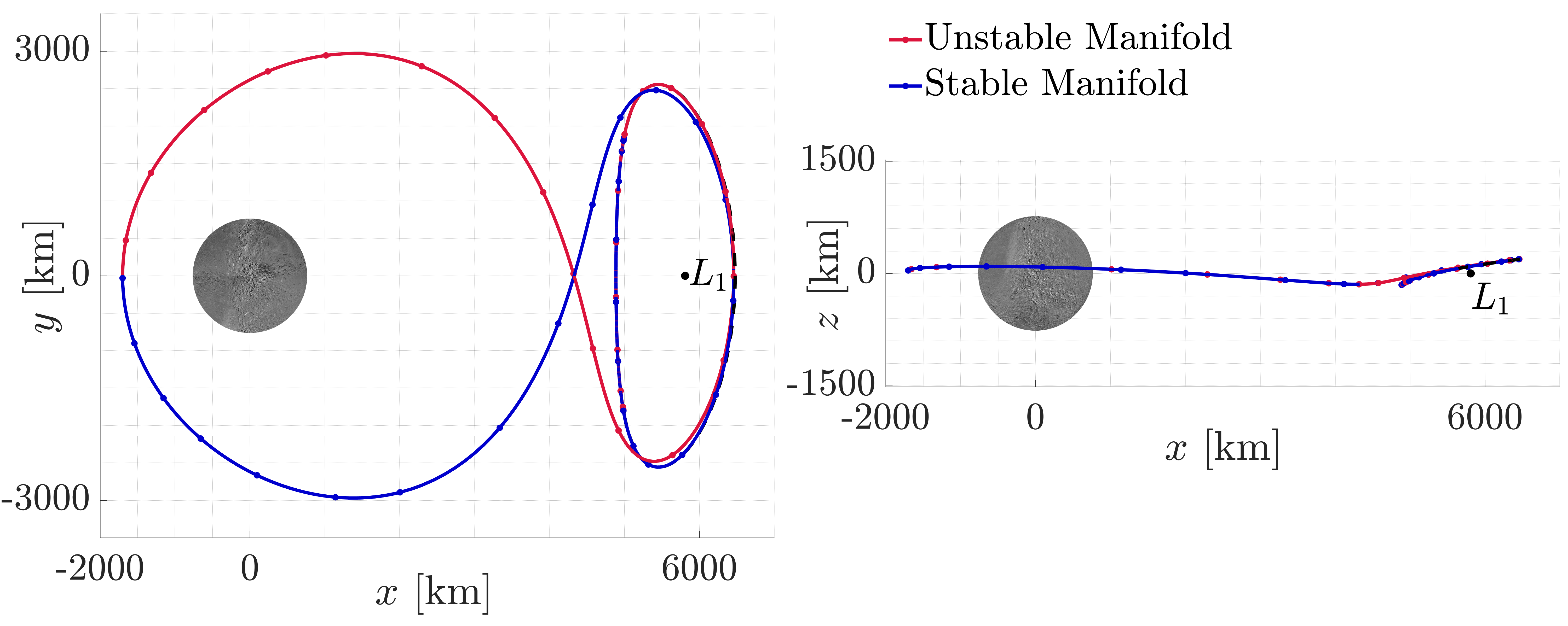}
       \caption{Planar projection of the representative homoclinic connection (Type E) for the Saturn–Rhea system. }
    \label{fig:connections12}
\end{figure*}

\begin{figure}[htpb]
    \centering
    \begin{subfigure}{\linewidth}
        \centering
        \includegraphics[width=0.58\linewidth]{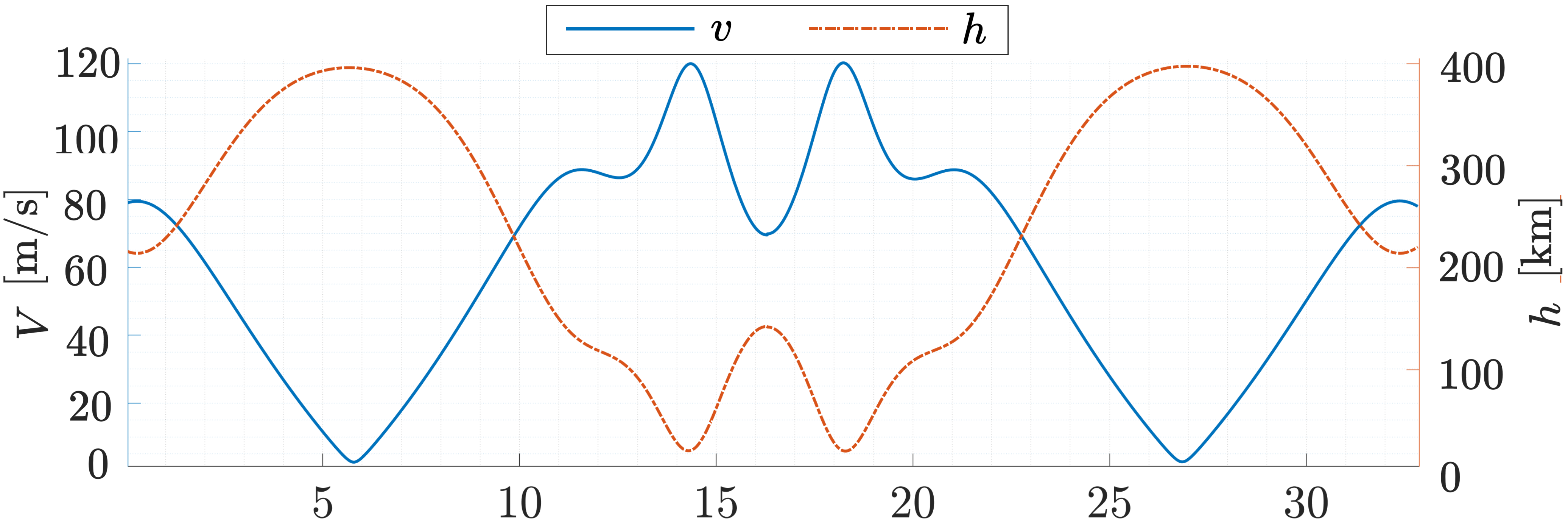}
    \end{subfigure}

    \begin{subfigure}{\linewidth}
        \centering
    \includegraphics[width=0.6\linewidth]{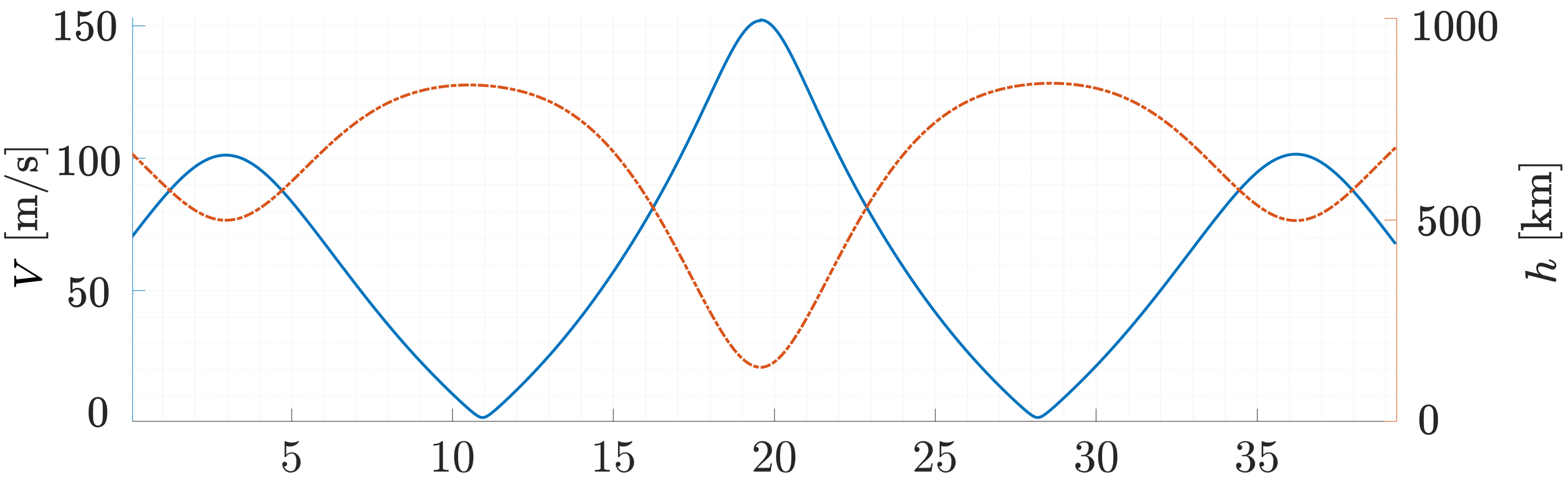}
    \end{subfigure}

    \begin{subfigure}{\linewidth}
        \centering
    \includegraphics[width=0.6\linewidth]{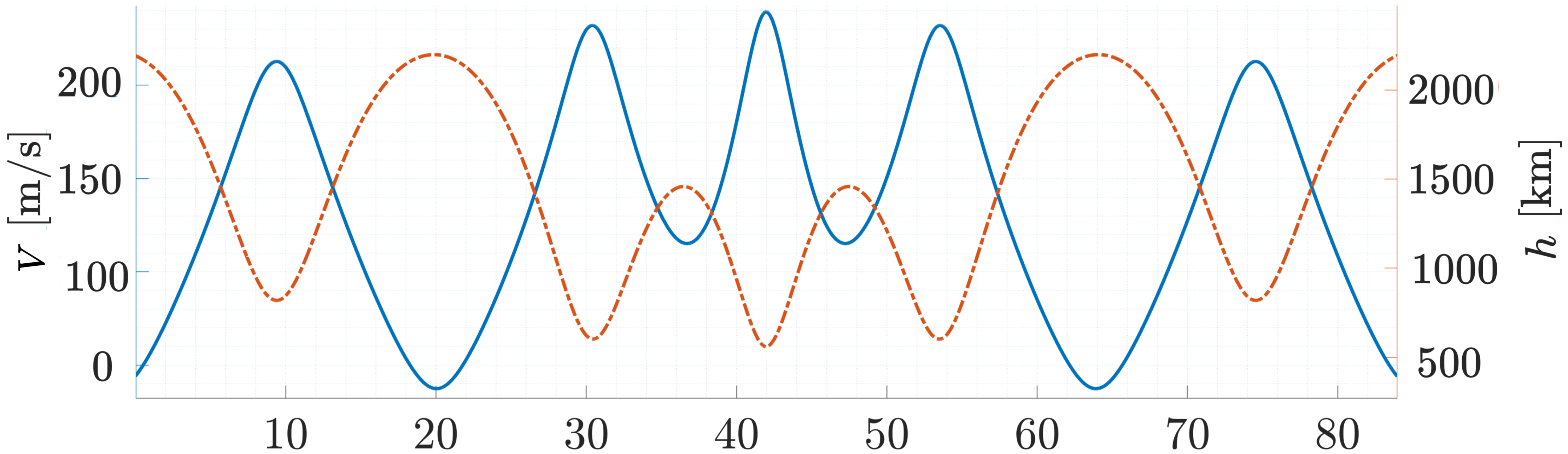}
    \end{subfigure}

    \begin{subfigure}{\linewidth}
        \centering
    \includegraphics[width=0.6\linewidth]{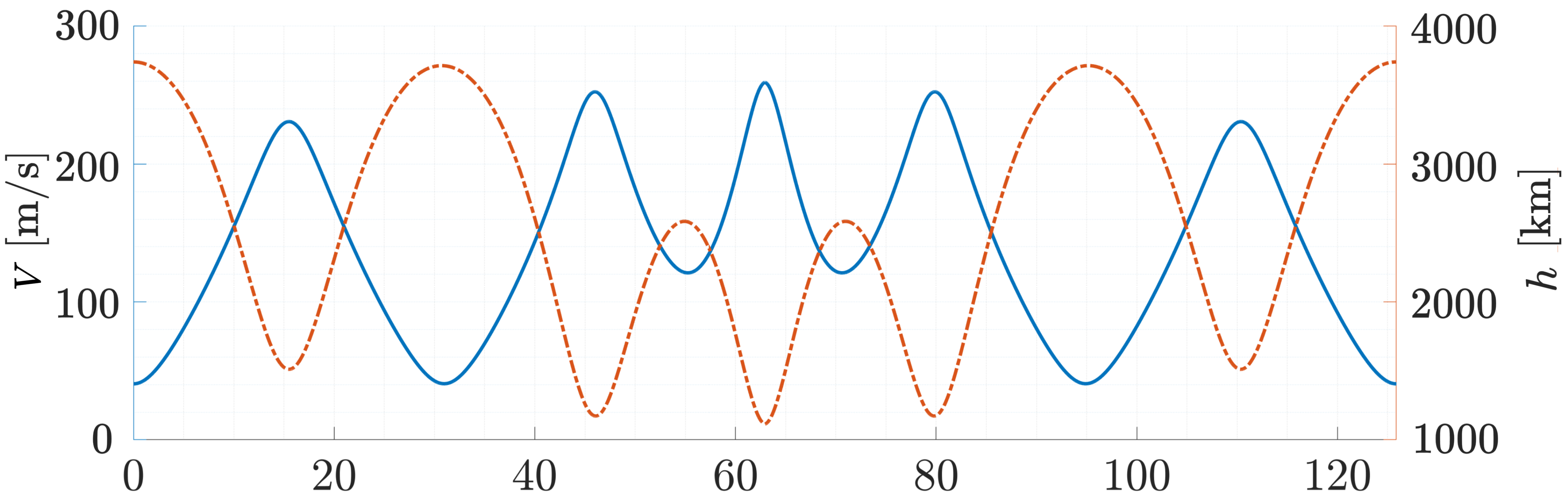}
    \end{subfigure}

    \begin{subfigure}{\linewidth}
        \centering
    \includegraphics[width=0.6\linewidth]{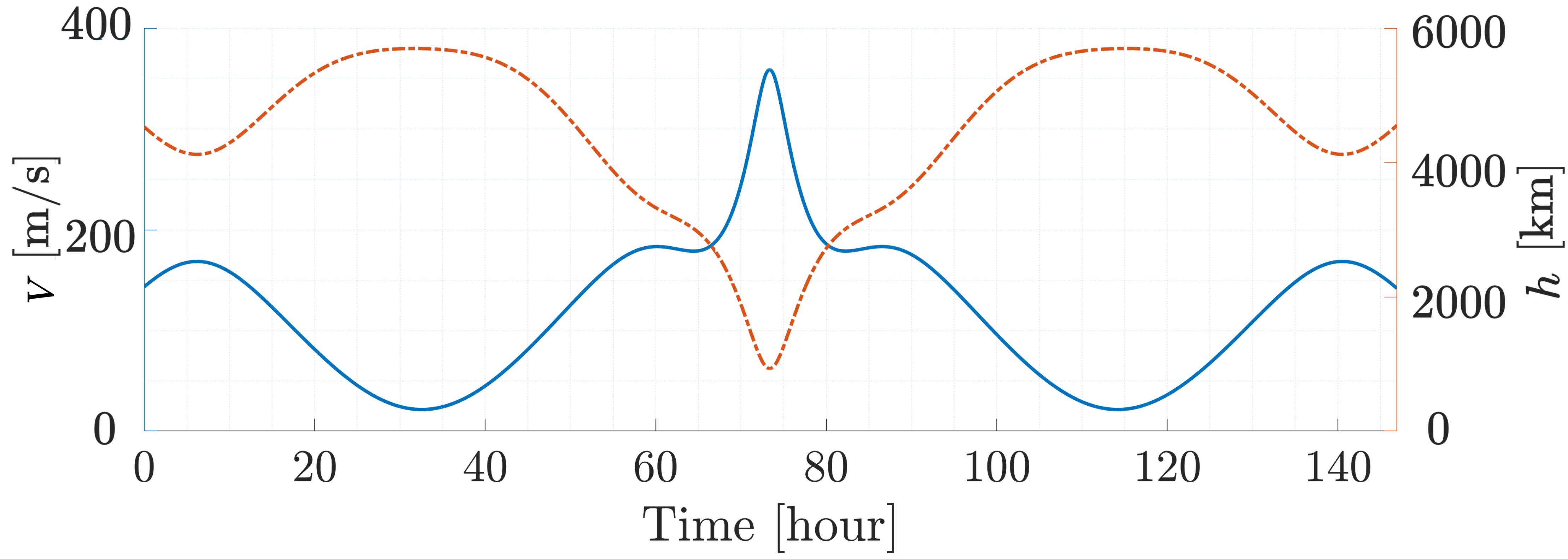}
    \end{subfigure}
     \caption{Time history of altitude and inertial velocity relative to the moon for the selected connections. From top to bottom: SMi (Type A), SEn (Type B), STe (Type C), SDi (Type D), and SRhe (Type E) systems.}
 \label{fig:hv}
\end{figure}

\begin{figure}[h]
    \centering
    \includegraphics[width=0.5\linewidth]{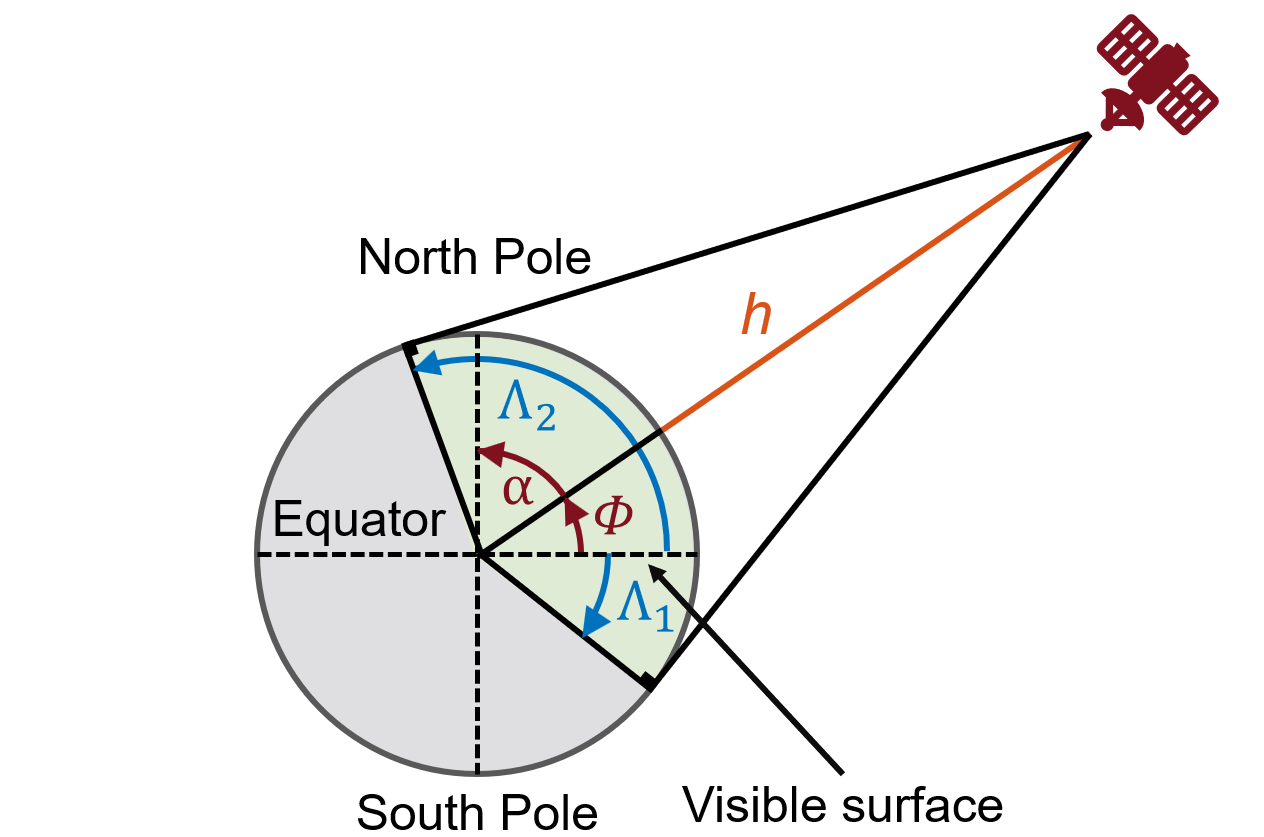}
    \caption{Instantaneous coverage parameters, $\Lambda_1$ and $\Lambda_2$. The shaded region (in green) indicates the portion of the surface instantaneously visible from the spacecraft.}
    \label{LambdaScheme}
\end{figure}

\begin{figure*}[htbp]
    \centering
    \begin{subfigure}{0.47\textwidth}
        \centering
        \includegraphics[width=\linewidth]{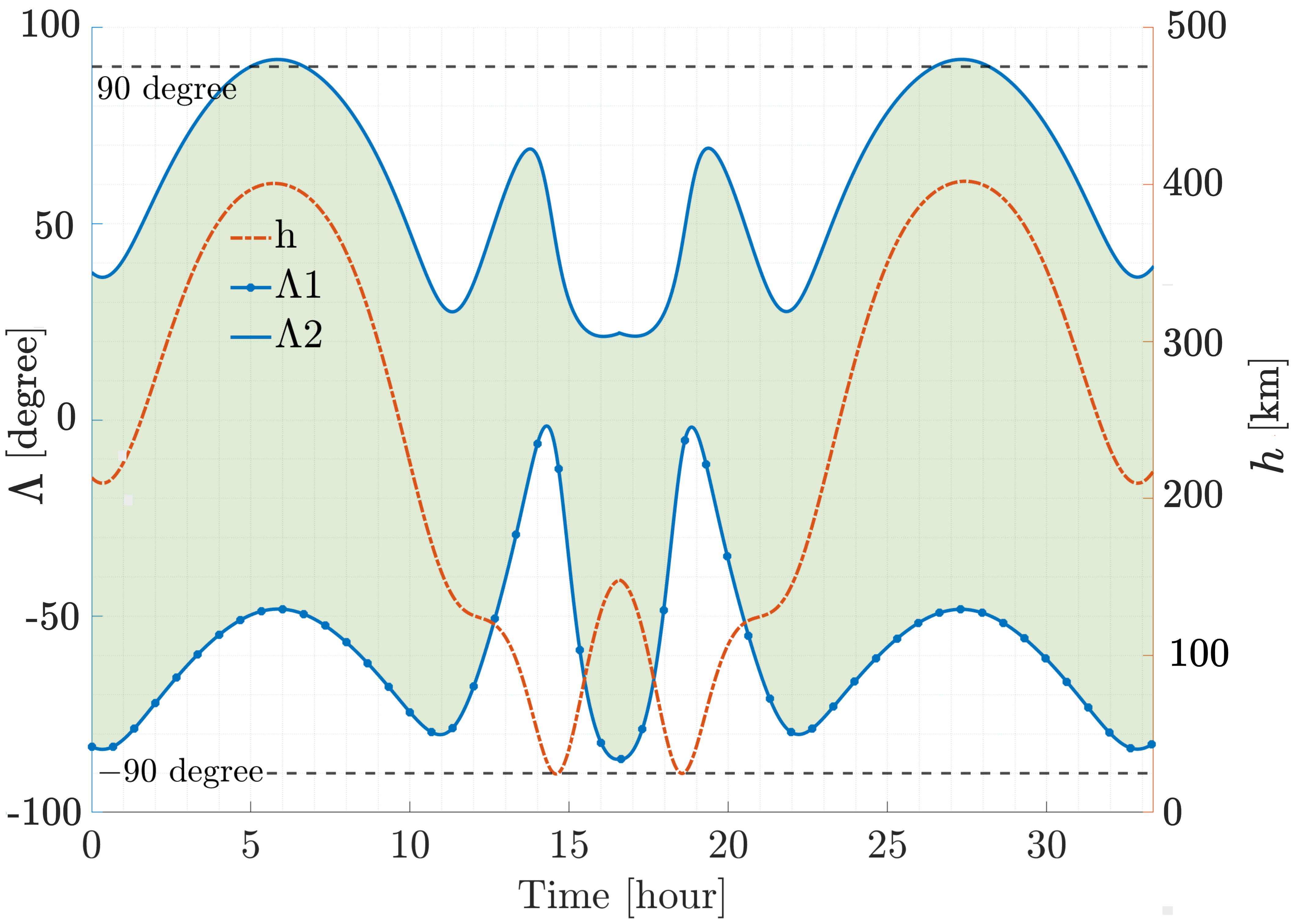}
         \caption{Type A - Saturn-Mimas system.}
    \end{subfigure}
    \begin{subfigure}{0.49\textwidth}
        \centering
        \includegraphics[width=\linewidth]{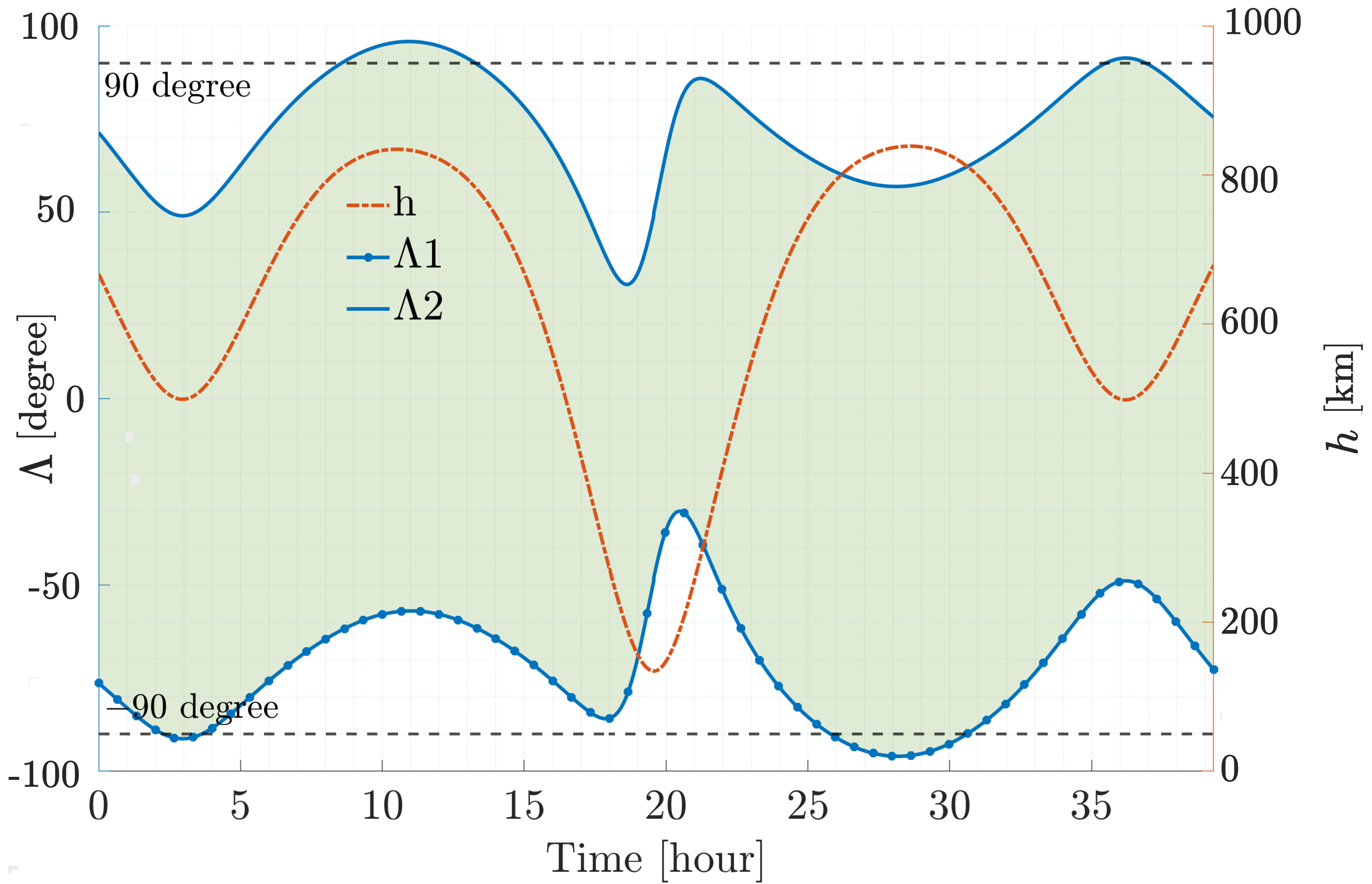}
        \caption{Type B - Saturn-Enceladus system.}
    \end{subfigure}
    \vspace{0.7cm}
    \begin{subfigure}{0.5\textwidth}
        \centering
        \includegraphics[width=\linewidth]{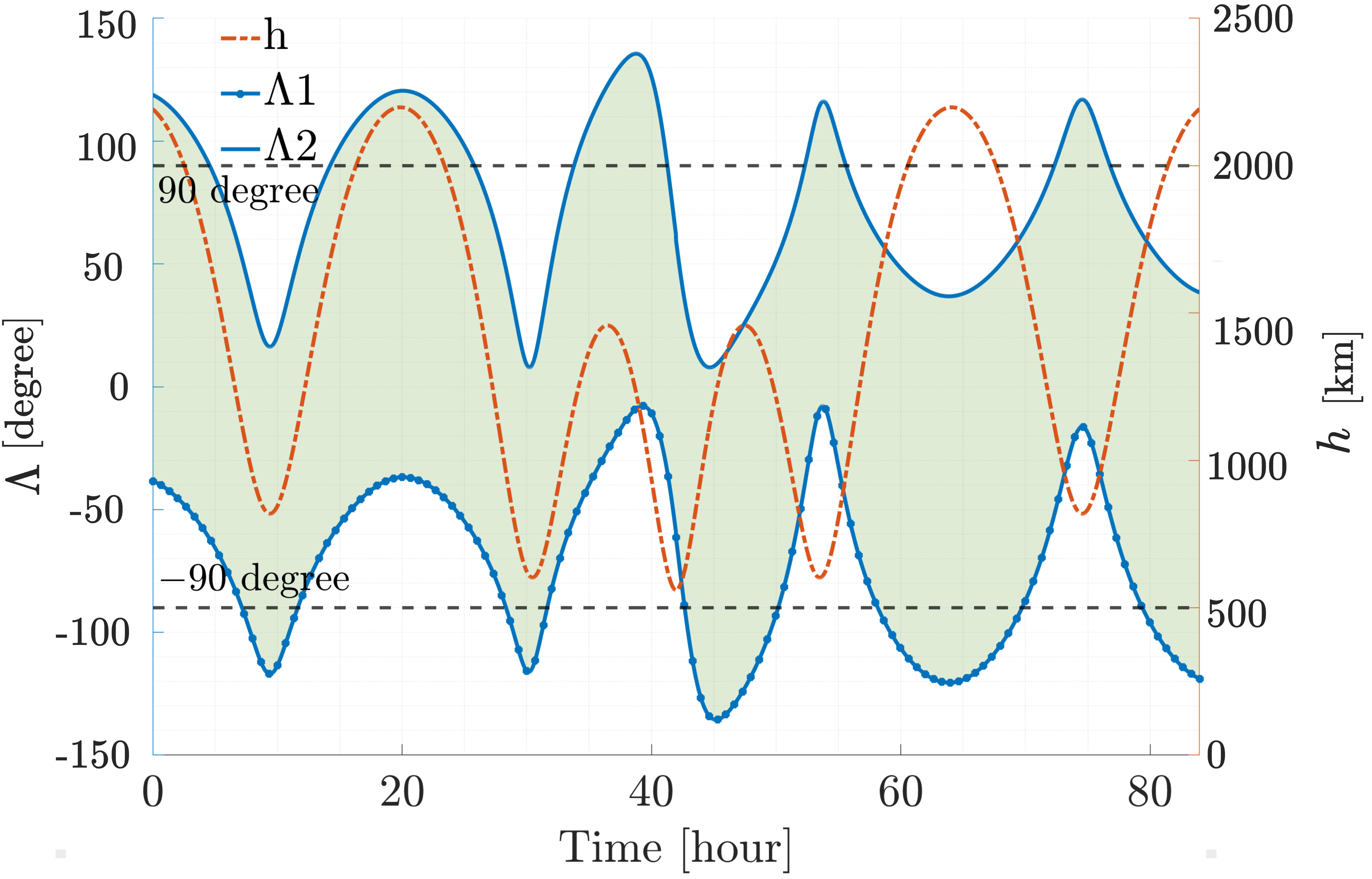}
        \caption{Type C - Saturn-Tethys system.}
    \end{subfigure}
    \begin{subfigure}{0.48\textwidth}
        \centering
        \includegraphics[width=\linewidth]{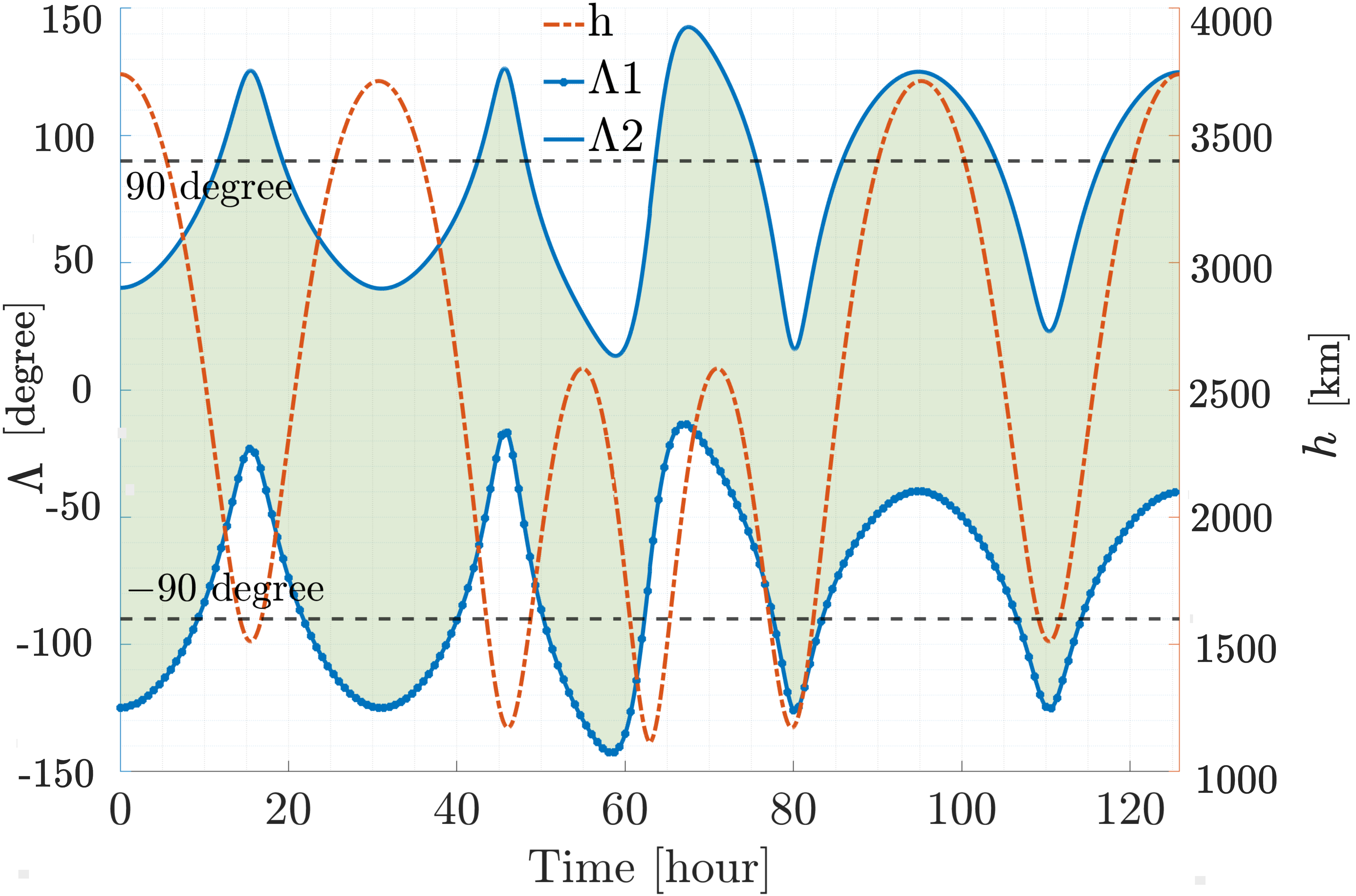}        
        \caption{Type D - Saturn-Dione system.}
    \end{subfigure}
    \begin{subfigure}{0.49\textwidth}
        \centering
        \includegraphics[width=\linewidth]{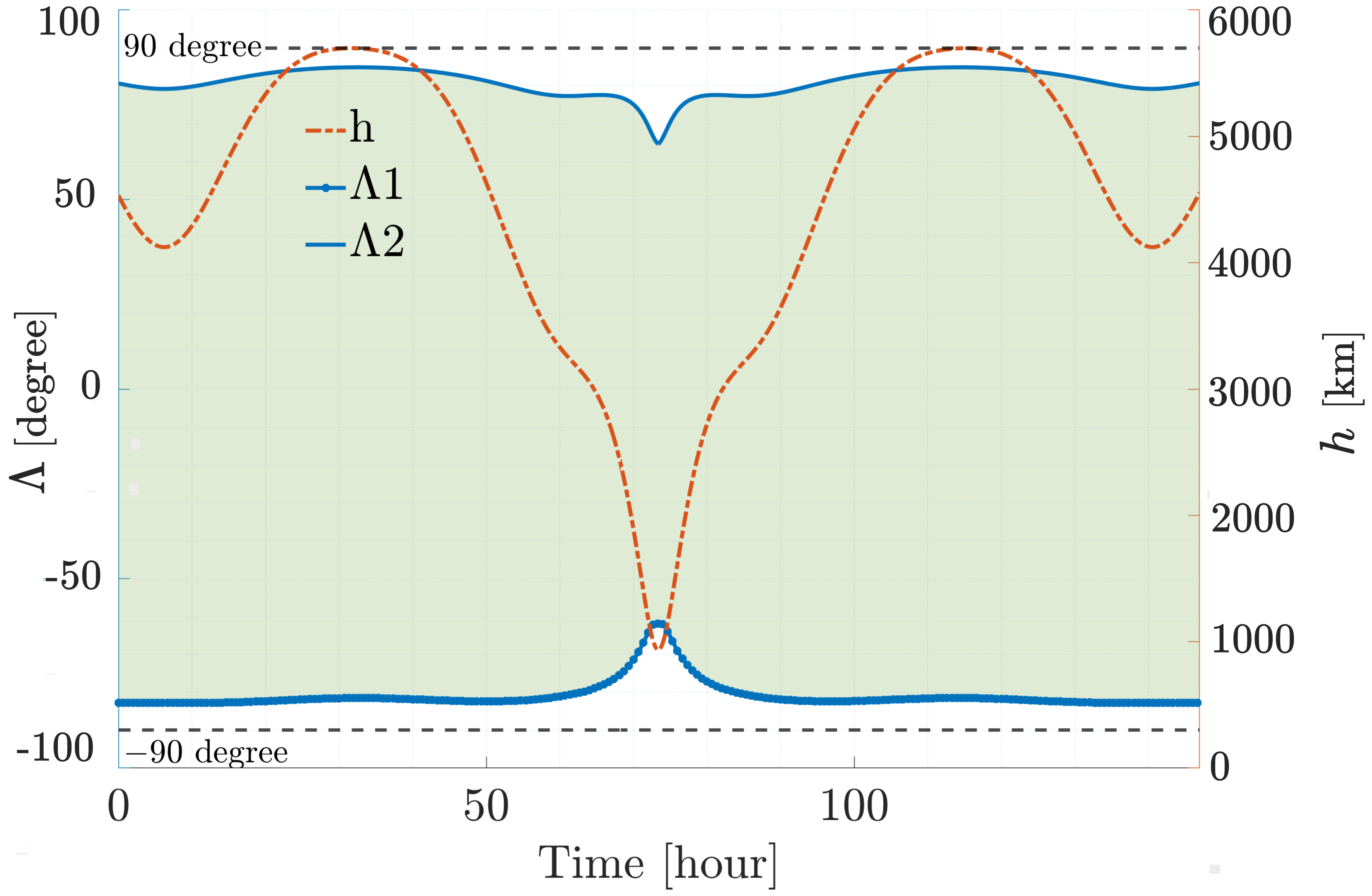}
        \caption{Type E -Saturn-Rhea system.}
    \end{subfigure}
    \caption{Time histories of $\Lambda_1$, $\Lambda_2$, and $h$ for the representative connections of Types A to E in the $J_2$-perturbed CR3BPs. Types A to D correspond to heteroclinic connections, while Type E is a homoclinic connection.}
    \label{fig:lambda}
\end{figure*}

\begin{figure}[H]
    \centering
    \begin{subfigure}{\linewidth}
        \centering
        \includegraphics[width=0.7\linewidth]{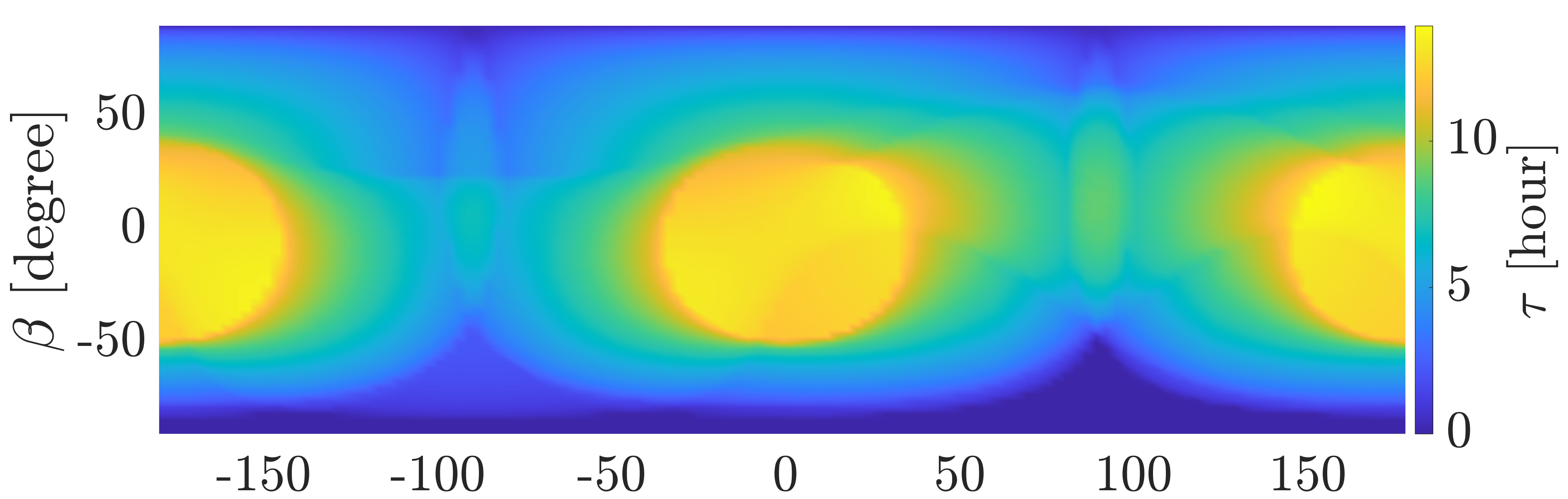}
    \end{subfigure}

    \begin{subfigure}{\linewidth}
        \centering
        \includegraphics[width=0.7\linewidth]{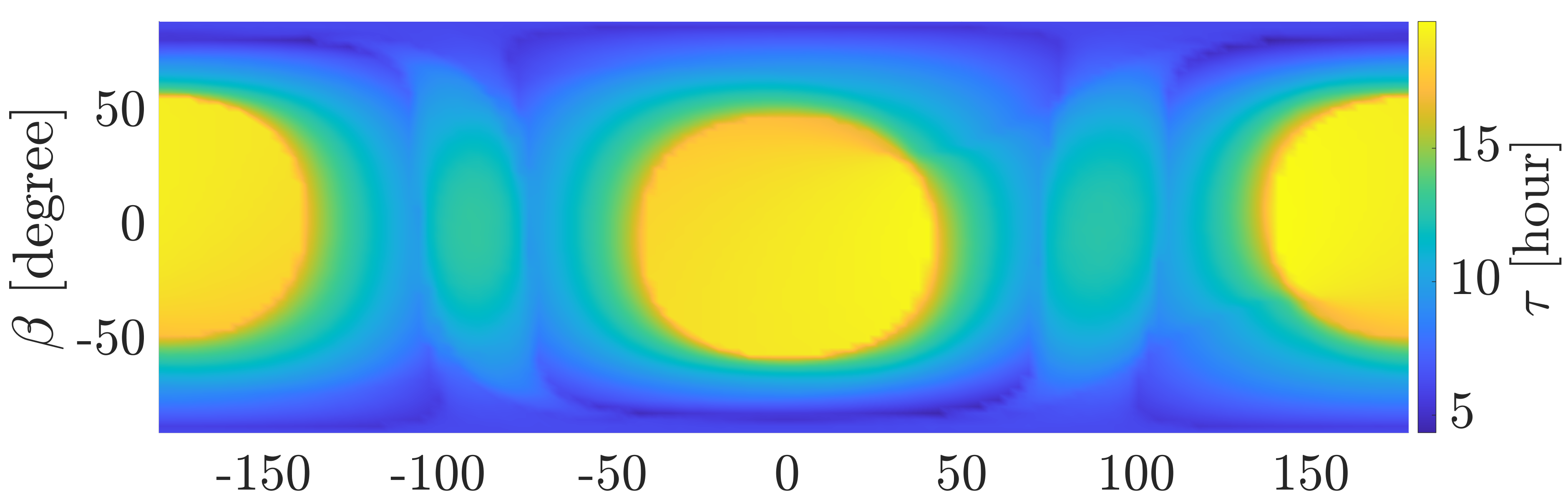}
    \end{subfigure}

    \begin{subfigure}{\linewidth}
        \centering
        \includegraphics[width=0.7\linewidth]{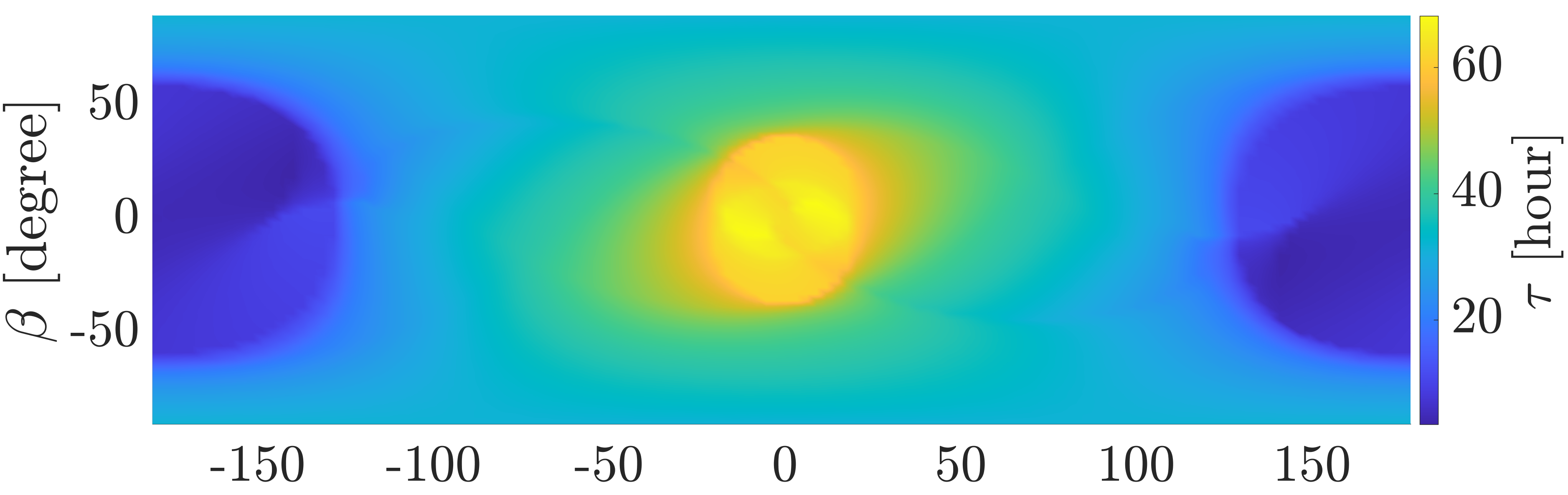}
    \end{subfigure}

    \begin{subfigure}{\linewidth}
        \centering
        \includegraphics[width=0.7\linewidth]{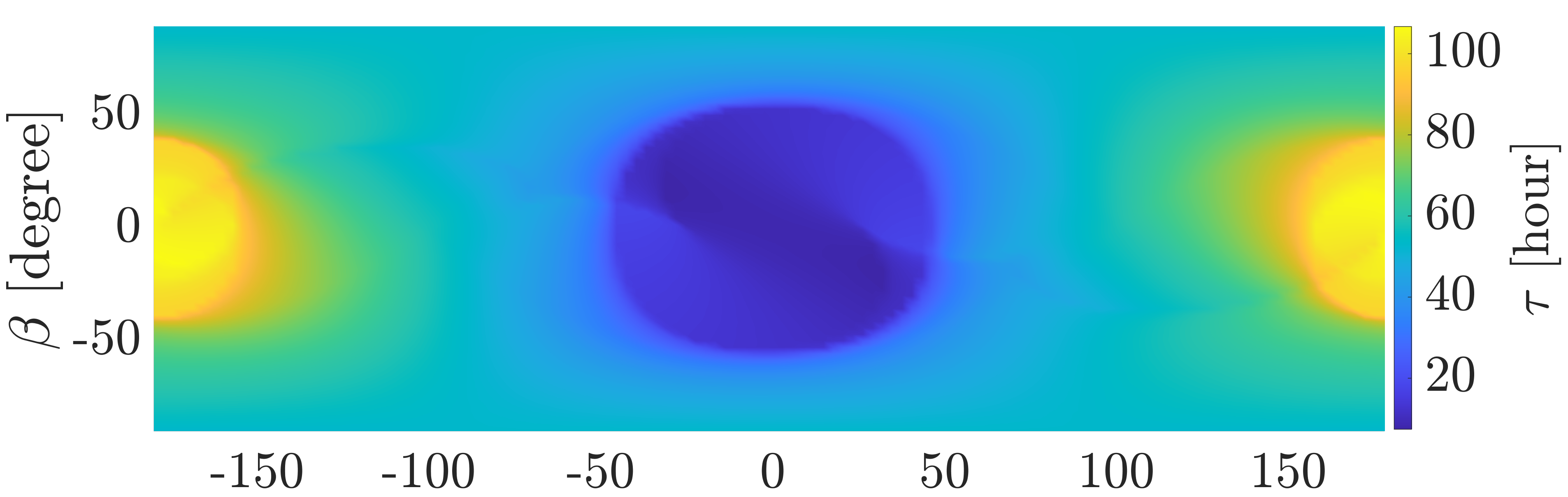}
    \end{subfigure}

    \begin{subfigure}{\linewidth}
        \centering
        \includegraphics[width=0.7\linewidth]{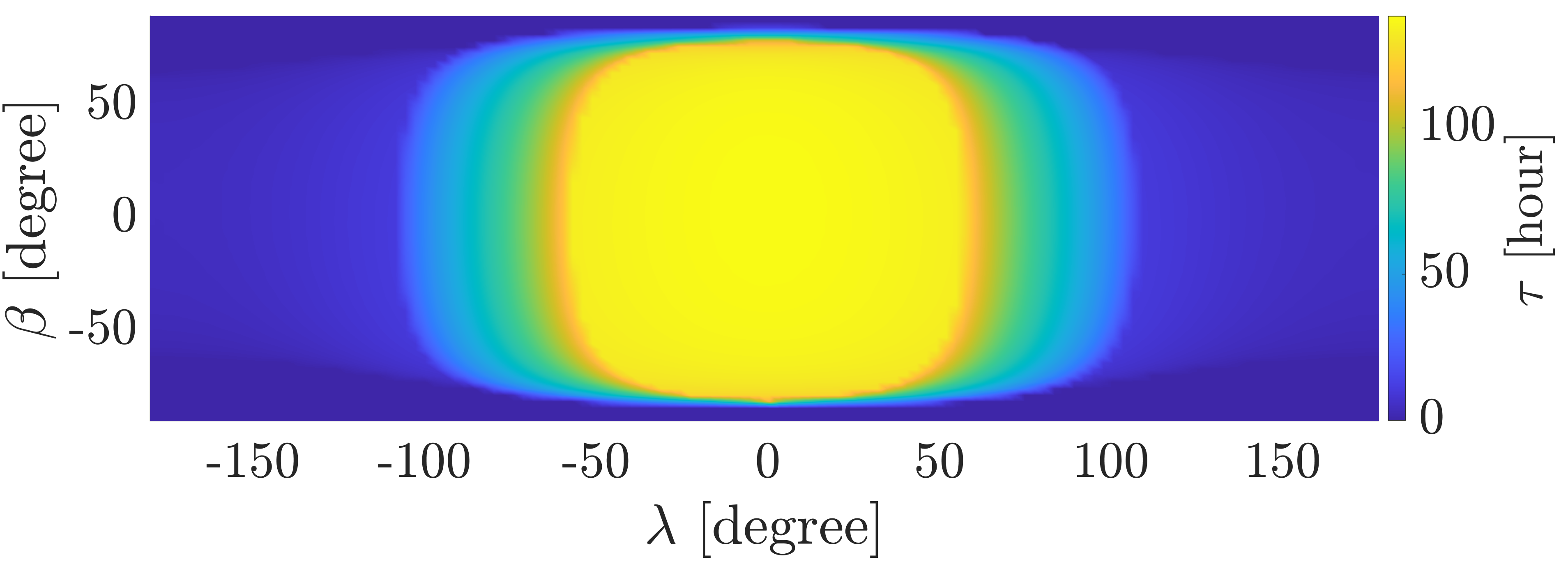}
    \end{subfigure}
    \caption{Geographical maps of total overflight time $\tau$ for the five selected connections. From top to bottom: SMi (Type A), SEn (Type B), STe (Type C), SDi (Type D), and SRh (Type E) systems.}
    \label{fig:tau}
\end{figure}

\noindent Each halo is discretized into 100 points, each point serving as the origin of an HIM trajectory. Homoclinic and heteroclinic connections are obtained with the method described in Section \ref{sciencedesign}. All connection types, including their corresponding inverse transfers, are successfully identified within the $J_2$-perturbed dynamical models for each ILM, except for Type B, E, and F connections in the Saturn–Mimas system. In these cases, the position and velocity discontinuities at $\Sigma$ significantly exceed the state-error thresholds adopted in this study (see Section \ref{sciencedesign}), and the corresponding trajectories are therefore excluded from the suitable set of connections.
The results emphasize the significant effect of the oblateness, particularly in limiting the number of connections with velocity errors below 1 m/s. The oblateness of the primary bodies introduces asymmetries in the manifolds, complicating the identification of viable low-energy transfers. Nevertheless, all identified connections exhibit position errors below 1 km and velocity errors less than 15 m/s, indicating near-continuity in both position and velocity at the Poincaré section. Strictly speaking, these transfers cannot be considered fully fuel-free; however, the magnitude of the mismatch implies that only a modest propellant consumption would be required.\\
\indent Only the most cost-efficient connections are reported. In particular, Table~\ref{tabresults_connections1} presents the results and characteristics for the five systems across the four types of heteroclinic connections (Types~A–D),
whereas Table~\ref{tabresults_connections2} summarizes the results for the two homoclinic types. In both tables, $\Delta T$ denotes the transfer time from the departure to the arrival halo (one full loop), $\Delta p$ the position error at~$\Sigma$, and $\Delta v$ the corresponding velocity error. For brevity, one representative connection type is illustrated for each Saturn–moon system: Type A for SMi, Type B for SEn, Type C for STe, Type D for SDi, and Type E for SRh, as shown in Figs. \ref{fig:connections1} and \ref{fig:connections12}. Type F connections are not reported, as their dynamical and geometrical characteristics are comparable to those of Type E, and would therefore not provide additional insights. Each trajectory is depicted in the corresponding moon-centered synodic frame. \\
\indent \indent The geometric and kinematic characteristics of an orbit determine its suitability for observational purposes. Key parameters include altitude, velocity, time of flight, instantaneous surface coverage, and length of visibility windows. Previous analysis of the observational performance of low-energy trajectories was conducted by Salazar et al. \cite{salazar2021science} for the Enceladus system only.
In this work, the study is extended to all five ILMs and to multiple types of connections, allowing a more comprehensive assessment of the observational capabilities of these trajectories. Figure \ref{fig:hv} illustrates the time history of the spacecraft altitude and the magnitude of its inertial speed \(V\) over the selected science orbits. \\
\indent To evaluate their observational performance, the instantaneous surface coverage of each moon is assessed using parameters \(\Lambda_1\) and \(\Lambda_2\), which define the limits of the central angle of coverage (\(2\alpha\)), measured positively northward from the equator. These parameters, illustrated in Fig.~\ref{LambdaScheme}, are given by\begin{align}
&\Lambda_1 = \phi - \alpha, \, \, \, \, \, \, \Lambda_2 = \phi + \alpha, \\
    &\text{with} \, \, \, \alpha = \cos^{-1} \left( \frac{R}{R + h} \right),
\end{align}
where \(h\) is the altitude and \(\phi\) the latitude of the spacecraft, measured from the equatorial plane of the moon and positive northward. Fig. \ref{fig:lambda} shows the time evolution of \(\Lambda_1\) and \(\Lambda_2\), with the green region highlighting the coverage amplitude, and horizontal dashed lines indicating the poles. \\
\indent The total overflight time \(\tau\), defined as the cumulative visibility of a specific surface point, is calculated for the entire moon surface along each transfer.  This parameter is determined by the visibility condition \(\theta_e \geq 0\), where \(\theta_e\) denotes the spacecraft elevation angle above the local horizon. To calculate \(\tau\), the flight time is divided into \(K\) segments of fixed duration \(\delta t\), with an elementary overflight time \(\delta \tau_i\) assigned to each segment \((i = 1, 2, \dots, K)\) as follows:
\begin{equation}
\delta \tau_i =
\begin{cases}
\delta t & \text{if } \theta_{e_i} \geq 0 \\
0 & \text{otherwise}
\end{cases} \longrightarrow \tau = \sum_{i=1}^{K} \delta \tau_i,
\end{equation}
\noindent where the total overflight time, for a given surface location, is obtained by summing the contributions from all segments. 
\begingroup
The computation of $\tau$ is performed in the moon-centered synodic reference frame, which rotates at the same rate as the moon about its spin axis, due to tidal locking with Saturn. As a consequence, the resulting  map quantifies the cumulative visibility time of each fixed surface point $(\lambda, \beta)$  over the duration of the orbit. 
\endgroup
The geographical maps of \(\tau\), presented in Fig.~\ref{fig:tau}, are generated by discretizing the surface of the moons with intervals of 1 degree in both longitude \(\lambda\) and latitude \(\beta\). \\
\indent For the science orbits considered in this work, the instantaneous coverage analysis based on $\Lambda_1$ and $\Lambda_2$ shows that, for connection types~B, C, and D, the entire surface of the corresponding moons is accessible. Different connection types exhibit markedly different polar performance. For example, as illustrated in Fig.~\ref{fig:lambda}, the Type B connection shown offers extended visibility of the polar regions, the most scientifically interesting areas, with up to 6 hours of coverage over the south pole of Enceladus, out of a total flight time of nearly 40 hours. For the Type~C and Type~D connections, the poles of Tethys and Dione remain visible for approximately half of the total flight time, distributed across windows at varying altitudes. In contrast, in the case of homoclinic connections (Type E), such as the one presented in Fig. \ref{fig:lambda}, and consistently observed across all Saturn–ILM systems, the surface coverage remains broad over time. However, the poles are not accessible because of the limited out-of-plane component of the trajectory, making them less suitable for detailed exploration of the Saturnian system. The corresponding overflight-time maps in Fig.~\ref{fig:tau} confirm that local visibility over the polar regions exceeds 6 hours for Types B, C, and D, while all connection types offer extended observation opportunities across broad equatorial bands.\\
\begingroup
\indent Quantitative coverage metrics for the representative science orbits are summarized in Table~\ref{tab:polar_coverage_metrics}. For each connection, the surface coverage is the fraction of the lunar surface for which $\tau> 0$. The polar-cap visibility metrics give, for the North (NP) and South (SP) poles, the total duration during which at least one point within the corresponding polar cap is in view. The maximum NP and SP revisit intervals represent the longest time gaps between successive visibility windows of the respective polar caps.\\
\endgroup
\begingroup
\indent It is worth emphasizing that, as already clear from Tables \ref{tabresults_connections1} and \ref{tabresults_connections2}, a broad variety of science-orbit geometries is available for each Saturn–ILM system, . In Figs.~\ref{fig:hv},\ref{fig:lambda} and\ref{fig:tau}, a single representative connection is shown for each moon, chosen to match the trajectories displayed in Fig.~\ref{fig:connections1} and \ref{fig:connections12}. However, by selecting different halo orbits and connection types, it is possible to obtain solutions with coverage metrics different from those reported in Table~\ref{tab:polar_coverage_metrics}, tailored to the desired scientific return. 
In particular, different choices of libration orbits allow extending the surface coverage, if required. This property holds for all the Moons, including Rhea, whose surface could be completely covered by choosing a different type of connection (e.g., Type B in place of Type E). The weak-capture strategy adopted in this work offers significant flexibility in the MT design and can enable full surface coverage by appropriately selecting the exploration orbits for each ILM.
\endgroup

\begin{table}[t]
    \centering
    \small
    \caption{{Coverage metrics for representative science orbits.}}
    \label{tab:polar_coverage_metrics}
    {\renewcommand{\arraystretch}{1.0}
    
    \begin{tabular}{lrrrrrr}
        \hline
        Connection &
        Surface &
        Time of &
        NP  &
        SP  &
        Max NP  &
        Max SP  \\
         &
        coverage &
        flight&
        visibility &
        visibility &
        revisit time&
        revisit time\\
         &
        [\%] &
        [hour] &
        [hour] &
        [hour] &
        [hour] &
        [hour] \\
        \hline
        Type A (SMi) &
        93.3 &
        32.83 &
        4.0 &
        0 &
        19.0 &
        -- \\
        Type B (SEn) &
        100.0 &
        39.32 &
        6.1 &
        6.2 &
        22.3 &
        22.2 \\
        Type C (STe) &
        100.0 &
        83.96 &
        32.8 &
        33.2 &
        16.4 &
        16.3 \\
        Type D (SDi) &
        100.0 &
        125.8 &
        58.3 &
        59.1 &
        23.6 &
        25.6 \\
        Type E (SRh) &
        82.4 &
        147.5 &
        0 &
        0 &
        -- &
        -- \\
        \hline
    \end{tabular}}
\end{table}

\subsection{Inter-moon transfers}\label{intmoonresults}
The inter-moon transfers leverage LT propulsion to patch the most favorable conditions among the stable and unstable HIMs of the periodic orbits used for moon exploration. For each transfer, 100 (points per halo) × 25 (halos per family) x 2 (families) conditions are available at each end, resulting in \( 2.5 \cdot 10^7\) candidate trajectories per segment. Each inter-SOI segment assumes that the $x$-axis of the synodic frame associated with the Saturn-departing moon system is aligned with the $X$-axis of the SCI frame for the computation of the departing osculating orbital elements. The initial spacecraft mass is set to 1000 kg at Rhea.  The tour starts at epoch 2042-January-01, consistent with the end of the interplanetary phase analyzed in~\cite{fantino2023end}, which brings the spacecraft into the Saturn system after a 12.34-year transfer (21~January~2028 to 24~May~2040) and a subsequent unpowered capture at Saturn via Titan flybys. All relevant perturbations described in Section~\ref{model} are accounted for in the propagation of each propelled arc, to enhance the fidelity of the results. In particular, the dynamical model includes the oblateness of Saturn, Titan and ILMs third-body gravitational perturbation, and the solar gravitational perturbation, only for the outermost transfers (Rhea–Dione and Dione–Tethys).

\begingroup
The Q-law in Eq.~\eqref{errfun} is implemented for the four LT arcs. The
elements of the weight vector $\boldsymbol{w} = [w_1, w_2, w_3, w_4]^T$ are initially set to unitary values and then refined by numerical tuning. In particular, the weights are optimized using the Nelder--Mead algorithm (\texttt{fminsearch} in
MATLAB), which efficiently explores the parameter space to identify values that satisfy the prescribed terminal tolerances while minimizing the transfer time (and thus, the propellant consumption). The segment-specific weight vectors are reported in Table~\ref{tab:Qlawweights} to facilitate reproducibility of the results.
\endgroup
\begin{table}[t]
    \centering
    \small
    \caption{Q-law weight vectors for the four inter-moon low-thrust arcs.}
    \label{tab:Qlawweights}
    \renewcommand{\arraystretch}{1.0}
    {
    \begin{tabular}{lc}
        \hline
        Transfer & $\boldsymbol{w}$ \\
        \hline
        Rh--Di & $[\,2.651,\; 0.442,\; 0.045,\; 0.356\,]$ \\
        Di--Te & $[\,2.458,\; 0.409,\; 0.042,\; 0.306\,]$ \\
        Te--En & $[\,2.403,\; 0.401,\; 0.041,\; 0.293\,]$ \\
        En--Mi & $[\,2.521,\; 0.421,\; 0.043,\; 0.322\,]$ \\
        \hline
    \end{tabular}}
\end{table}

Among all feasible solutions, the trajectory minimizing propellant usage, expressed through velocity variation, is selected for each transfer leg. The performance parameters of the optimal propelled arcs are summarized in Table~\ref{tab:inter_moon_transfers}. It is worth noting that the four transfers are executed sequentially, with the initial mass of each leg corresponding to the final mass of the previous one. 

\begin{table}[h!]
\caption{Magnitude of total impulse provided \(\Delta V\), time of flight \(\Delta t\), mass consumption \(\Delta m\) for the optimal LT inter-moon transfers between SOIs.}
\centering
\renewcommand{\arraystretch}{0.8}
\begin{tabular}{l r r r}
\hline
{Transfer} & \(\Delta V\) [m/s] & \(\Delta t\) [day] & \(\Delta m\) [kg] \\
\hline
Rh-Di & 1084  & 336 & 66.8 \\
Di-Te & 876  & 246 & 48.9 \\
Te-En & 958  & 254 & 50.3 \\
En-Mi & 1438 & 331 & 65.5 \\
\textbf{Total} & \textbf{4356} & \textbf{1167} & \textbf{231.5} \\
\hline
\end{tabular}
\label{tab:inter_moon_transfers}
\end{table}

\begin{figure*}[h]
    \centering
    \begin{subfigure}{0.4\textwidth}
        \centering
        \includegraphics[width=\linewidth]{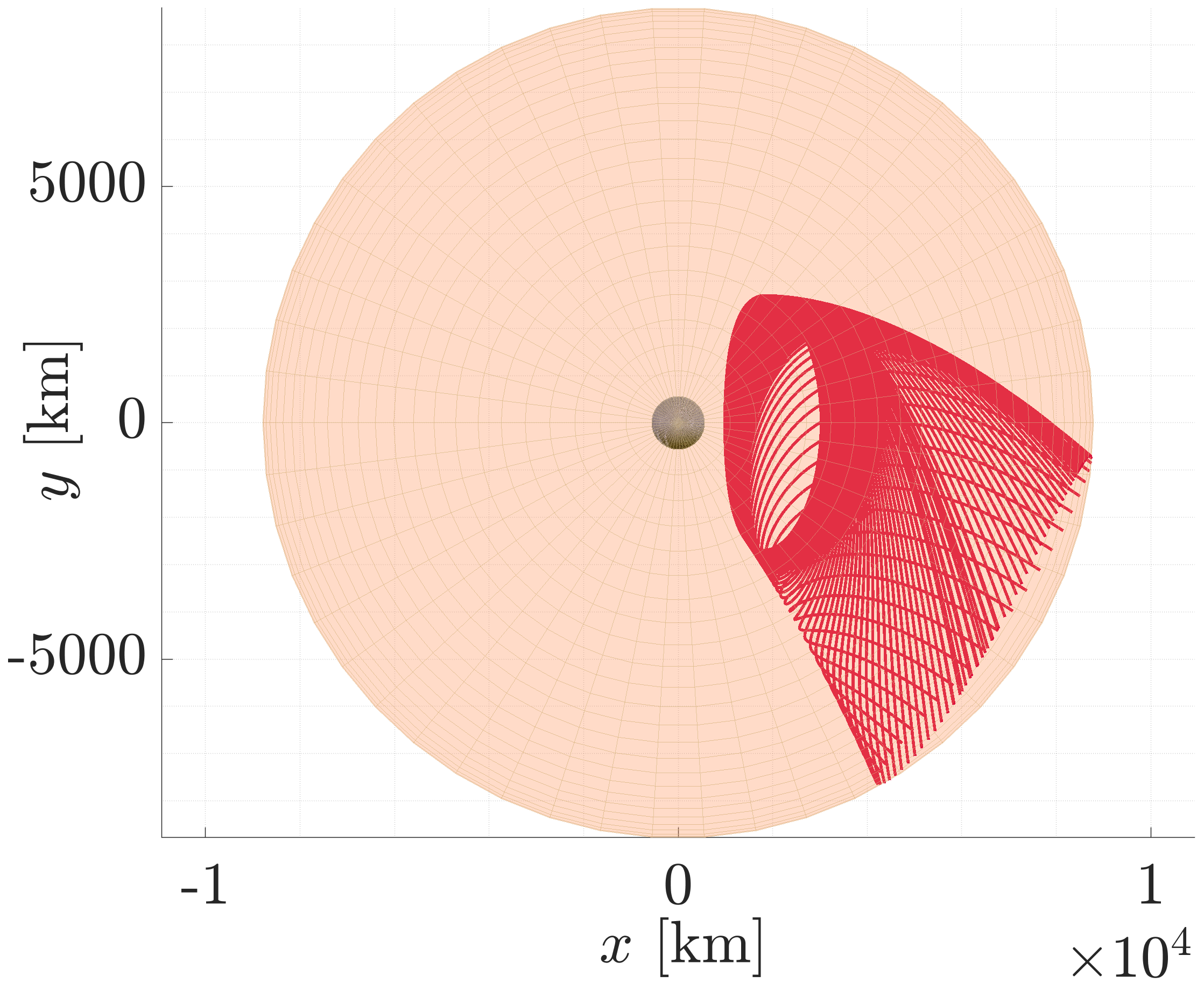}
    \end{subfigure}
\hspace{0.02\textwidth}
    \begin{subfigure}{0.4\textwidth}
        \centering
        \includegraphics[width=\linewidth]{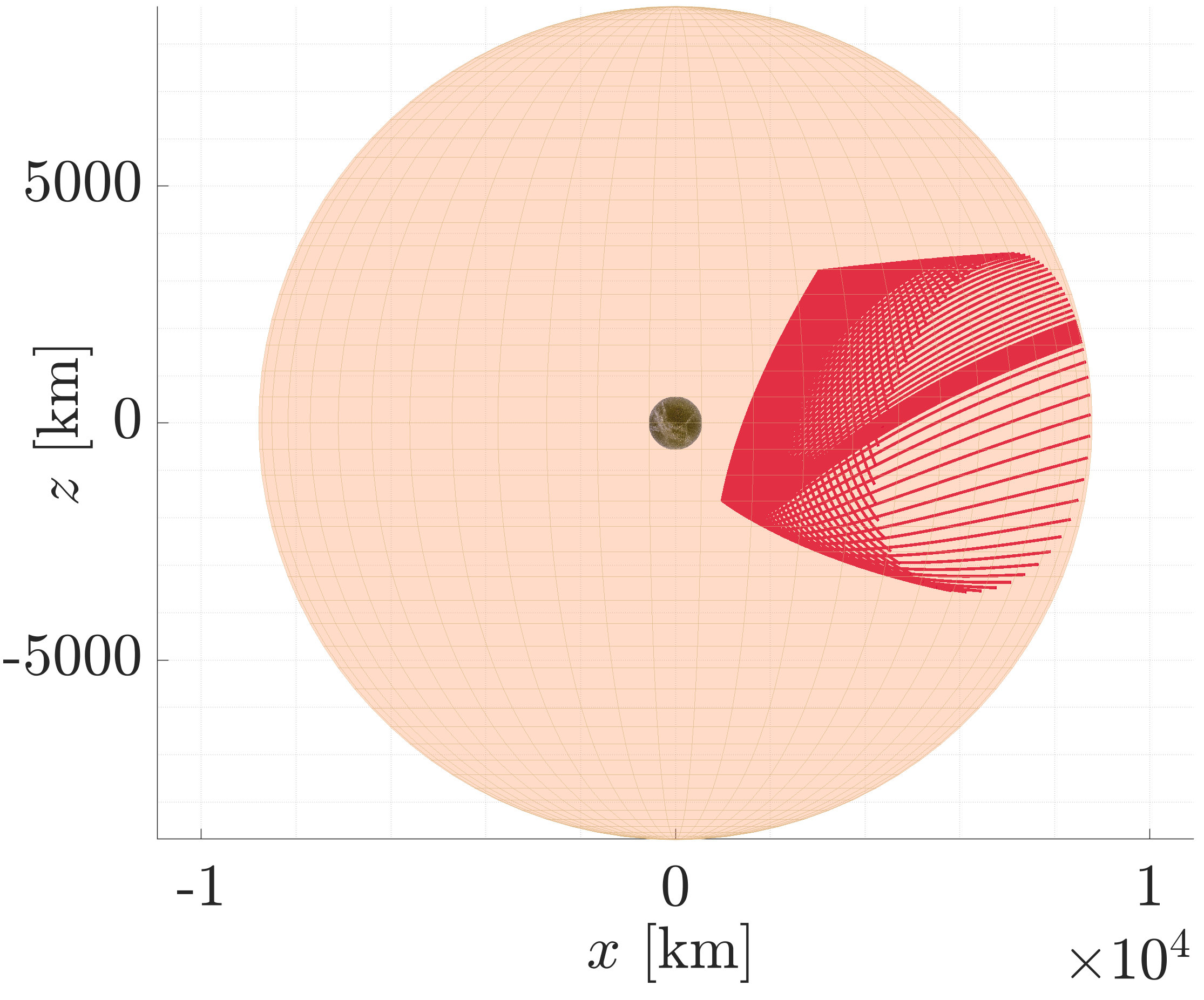}
    \end{subfigure}
    \caption{Global unstable manifold (outward branch) of a halo orbit around $L_1$ of the Saturn–Dione system and its intersections with the surface of the SOI. Left: $xy$ projection; right: $xz$ projection.}
    \label{fig:dione_out}
\end{figure*}

\begin{figure*}[h]
    \centering
    \begin{subfigure}{0.4\textwidth}
        \centering
        \includegraphics[width=\linewidth]{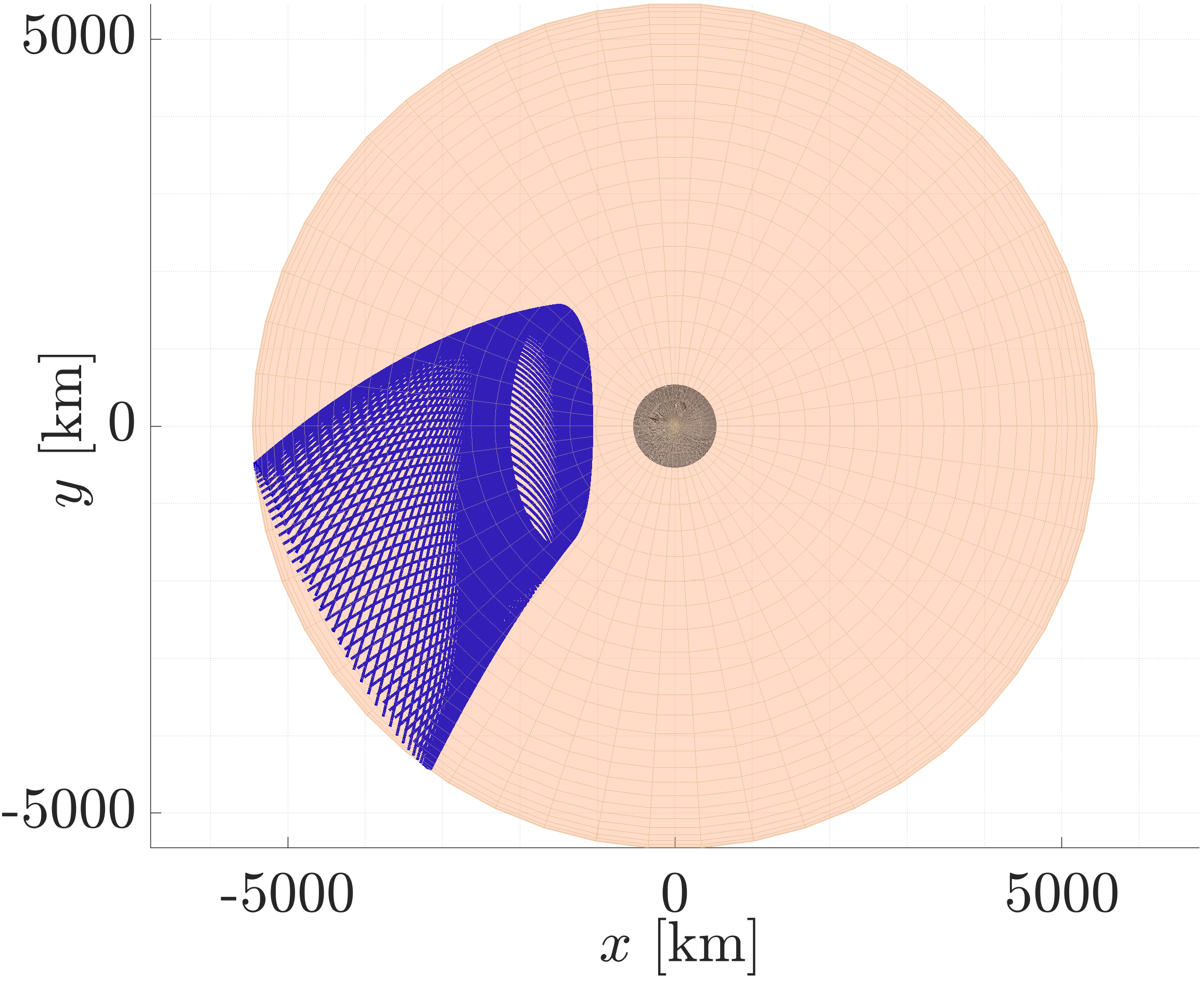}
    \end{subfigure}
\hspace{0.02\textwidth}
    \begin{subfigure}{0.4\textwidth}
        \centering
        \includegraphics[width=\linewidth]{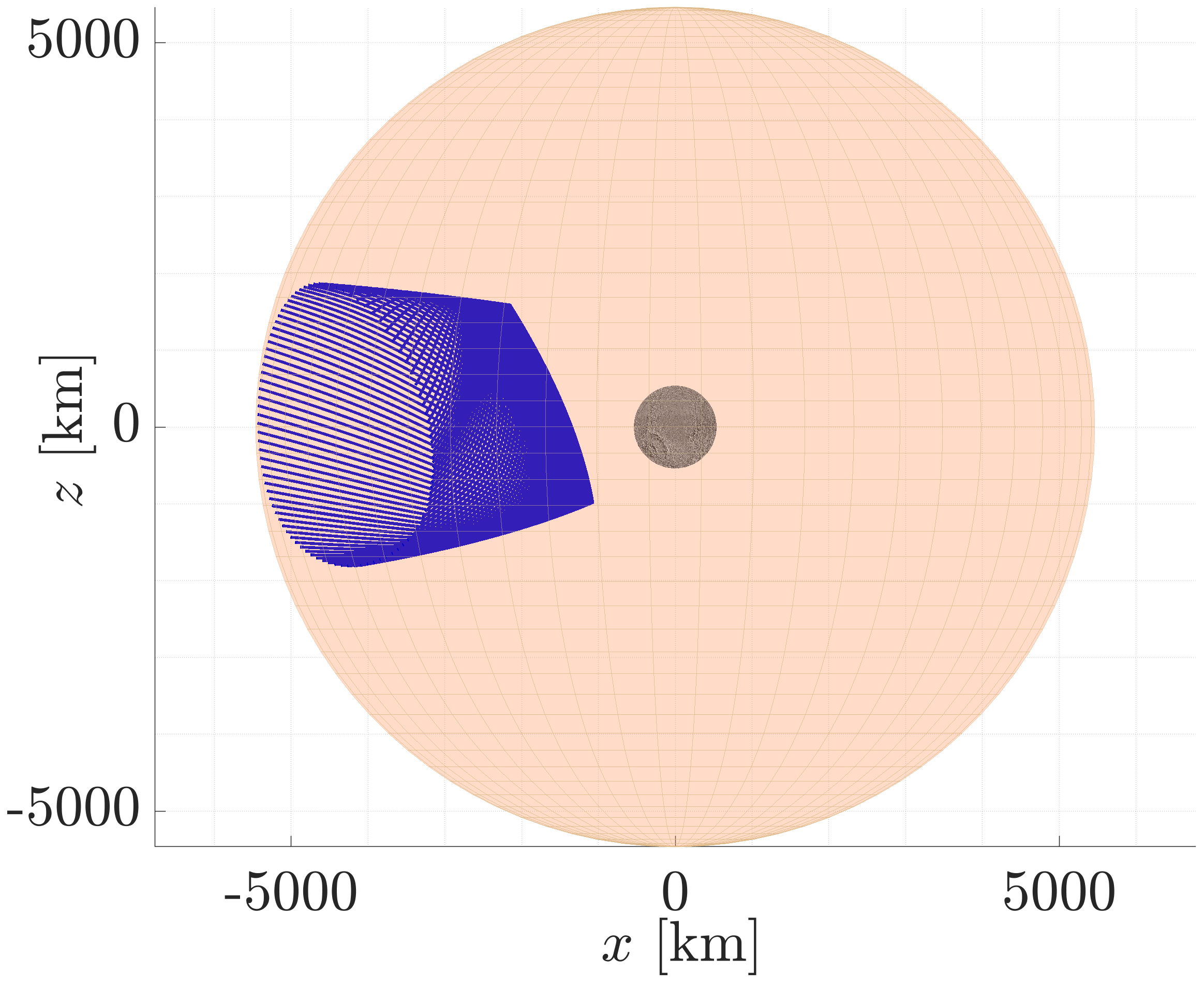}
    \end{subfigure}
\caption{Global stable manifold (inward branch) of a halo orbit around $L_2$ of the Saturn–Tethys system and its intersections with the surface of the SOI. Left: $xy$ projection; right: $xz$ projection.}
\label{fig:teti_in}
\end{figure*}

\indent Although most of the flight time is spent in the Saturnian inter-moon space, each transfer also includes a short-duration phase within the SOIs of the departure and arrival moons. These segments correspond to the manifold trajectories linking the halo orbits to the respective SOI boundaries. Table~\ref{tab:time2soi} summarizes the corresponding transfer durations, reporting the minimum and maximum flight times across the family (columns 4 and 5), the value associated with the optimal solution (i.e., the transfer with the lowest $\Delta V$ of the propelled arc, column 6), and the corresponding halo index and family type (columns 7 and 8). The results highlight the variability of the intra-SOI transfer duration to the specific halo orbit employed. Figures~\ref{fig:dione_out} and~\ref{fig:teti_in} show representative examples of the global stable and unstable manifolds within the moons’ respective SOIs for the Dione–Tethys transfer.\\
\begin{table*}[t]
\caption{Transfer time between departure or arrival halo orbit and the SOI for all the inter-moon transfers.}
\centering
\small
{\renewcommand{\arraystretch}{1.1}
\begin{tabular}{l l l r r r r r}
\hline
{System} & \(L_i\) & {Stability} & \(\Delta t^\text{min}\) & \(\Delta t^\text{max}\) & \(\Delta t^\text{opt}\) & {Optimal } & {Optimal} \\
{} & - & - & [day] & [day] & [day] & {orbit index} & {orbit family} \\
\hline
SRh & \(L_1\) & Unstable   & 3.36 & 22.23 & 3.56 & 23 & South\\
SDi & \(L_2\) & Stable   & 1.99 & 13.68 & 2.55 & 19 & South \\
SDi & \(L_1\) & Unstable & 1.98 & 13.52 & 2.08 & 24 & North \\
STe & \(L_2\) & Stable   & 1.33 & 8.16 & 1.51 & 20 & North \\
STe & \(L_1\) & Unstable & 1.34 & 6.88 & 4.95 & 1 & South  \\
SEn & \(L_2\) & Stable   & 0.88 & 1.59 & 1.02 & 19 & South \\
SEn & \(L_1\) & Unstable & 0.89 & 2.03 & 1.15 & 12 & North \\
SMi & \(L_2\) & Stable   & 0.57 & 1.06 & 1.04 & 2 & North\\
\hline
\end{tabular}}
\label{tab:time2soi}
\end{table*}
\indent Because the natural dynamics of the $J_2$-perturbed CR3BP are used within each SOI, the cost of the moon-to-moon transfer corresponds solely to the powered legs. The low-thrust portion of the MT requires a total \(\Delta V\) of 4356 [m/s], with a duration of 1167 days ($\approx$ 3.2 years) and a propellant consumption of 231.5 kg, resulting in a final spacecraft mass of 768.5~kg. During the four segments, the spacecraft performs approximately 100, 110, 160, and 290 revolutions around Saturn, offering opportunities to study the E ring and its surrounding environment.



\begingroup 

\subsection{End-to-end mission budget and timeline}

\begin{table*}[t]
    \centering
    \small
    \caption{Representative end-to-end MT phases with associated cost, duration, and mass evolution.}
    \label{tab:mission_budget}
    {\renewcommand{\arraystretch}{0.9}
    
    \begin{tabular}{l c r r r r r}
        \hline
        Phase & Start date & $\Delta t$ & $\Delta V$ & Thrust time & $\Delta m$ & Final mass \\
         & [yyyy mm dd] & [day] & [m/s] & [day] & [kg] & [kg] \\
        \hline
        SRh Type E   & 2042-01-01 & 6.13  & $\approx 0$ & $\approx 0$ & $\approx 0$ & 1000.0 \\
        Rh SOI exit     & 2042-01-07 & 3.56  & 0        & 0        & 0        & 1000.0 \\
        Rh--Di          & 2042-01-10 & 336   & 1084     & 336      & 66.8     & 933.2  \\
        Di SOI entry    & 2042-12-12 & 2.55  & 0        & 0        & 0        & 933.2  \\
        SDi Type D & 2042-12-15 & 5.25  & $\approx 0$ & $\approx 0$ & $\approx 0$ & 933.2  \\
        Di SOI exit     & 2042-12-20 & 2.08  & 0        & 0        & 0        & 933.2  \\
        Di--Te          & 2042-12-22 & 246   & 876      & 246      & 48.9     & 884.3  \\
        Te SOI entry    & 2043-08-25 & 1.51  & 0        & 0        & 0        & 884.3  \\
        STe Type C & 2043-08-27 & 3.50  & $\approx 0$ & $\approx 0$ & $\approx 0$ & 884.3  \\
        Te SOI exit     & 2043-08-30 & 4.95  & 0        & 0        & 0        & 884.3  \\
        Te--En          & 2043-09-04 & 254   & 958      & 254      & 50.3     & 834.0  \\
        En SOI entry    & 2044-05-15 & 1.02  & 0        & 0        & 0        & 834.0  \\
        SEn Type B & 2044-05-16 & 1.64  & $\approx 0$ & $\approx 0$ & $\approx 0$ & 834.0  \\
        En SOI exit     & 2044-05-18 & 1.15  & 0        & 0        & 0        & 834.0  \\
        En--Mi          & 2044-05-19 & 331   & 1438     & 331      & 65.5     & 768.5  \\
        Mi SOI entry    & 2045-04-15 & 1.04  & 0        & 0        & 0        & 768.5  \\
        SMi Type A & 2045-04-16 & 1.36  & $\approx 0$ & $\approx 0$ & $\approx 0$ & 768.5  \\
        \textbf{Total}  & \textbf{2045-04-18}         & \textbf{1203} & \textbf{4356} & \textbf{1167} & \textbf{231.5} & \textbf{768.5} \\
        \hline
    \end{tabular}}
\end{table*}

The tour duration reported in Section \ref{intmoonresults} refers exclusively to the cumulative $\Delta t$ of the powered transfers. To obtain the duration of the complete MT, the flight times associated with the halo-to-SOI legs and with the science orbits reported in Section \ref{resultsSO} must be added. Table~\ref{tab:mission_budget} reports the end-to-end mission budget and timeline for a representative realization of the proposed tour from the science orbit in the Saturn-Rhea system to the one around Mimas. The table includes the dates, $\Delta V$, thrust duration, and spacecraft mass evolution. Propulsive costs are associated exclusively with the low-thrust inter-moon transfers. Homoclinic and heteroclinic connections are assumed to have a negligible impact on the overall mission budget, thanks to the use of natural dynamics and the marginal $\Delta V$ required. Moreover, the halo-to-SOI phases are fully ballistic, as they follow the natural dynamics of the $J_2$-perturbed CR3BP. In defining the tour timeline, it is assumed that the spacecraft reaches the arrival halo orbit and departs immediately along the heteroclinic/homoclinic connection leading to the departure halo of the next inter-moon segment. Under this assumption, the tour starts on 2042-01-01 and ends on 2045-04-18, 
with a total duration of 1203 days ($\approx$ 3.3 years), and a final mass 
of 768.5~kg. For the adopted propulsion system (see Table~\ref{tab:propulsion}), the cumulative
thrust duration of 1167~days corresponds to a continuous power demand of 640~W during the low-thrust arcs, which is compatible with existing RTGs-powered electric-propulsion technologies for deep-space missions.

In some cases, the halo orbits providing the most favorable gateways for the moon-to-moon transfers do not coincide with those selected for the observational phases. When this occurs, suitable transfers between halo orbits of the same family are required. These transfers have marginal propellant cost and short flight times, and therefore barely affect the overall mission $\Delta V$ budget and total duration \cite{hiday1996impulsive,gomez1998study}. In this context, phasing maneuvers between successive connections are also neglected, as they do not affect the values reported in Table~\ref{tab:mission_budget}.

One of the key advantages of leveraging low-energy transfers for the intra-moon phases is the ability to modulate the mission timeline, allowing repeated science orbits with only marginal propellant consumption. For instance, if each science orbit were repeated 100 times, the total mission duration would increase to approximately 5 years. Such flexibility allows the mission timeline to be tailored to the desired scientific return, enabling multi-epoch observations of the moons and the long-term monitoring of dynamic phenomena such as subsurface activity or Enceladus’ geyser-like jets. This degree of adaptability is not achievable with traditional impulsive or high-energy transfer strategies. The MT outlined here should therefore be interpreted as an illustrative end-to-end example enabled by the proposed strategy, with the very interesting prospect of adjusting and refining the phases in relation to contingent mission specifications and objectives.
\endgroup

\section{Comparative analysis and discussion of results}\label{disc}

This work advances a tour-design strategy based on the combined use of low-energy dynamics and low-thrust propulsion within a perturbed dynamical environment. It is a very significant extension and refinement of the mission concept introduced by Fantino et al. in~\cite{fantino2023end}. In particular, the present study improves upon the unperturbed 2D tour proposed, in which the science orbits were designed in the planar CR3BP and the moon-to-moon transfers were computed in a 2BP framework. 
Here, the perturbations induced by the oblateness of Saturn and of each target ILM are included in the intra-SOI CR3BP dynamics, while the inter-SOI legs are propagated in an $N$-body model that also incorporates Saturn’s $J_2$. The complete tour is entirely three-dimensional, enabling full-surface coverage during the science-orbit phases. The proposed 3D tour achieves performance comparable to the planar 2D solution (characterized by a total $\Delta V$ of 3.1~km/s for three moon-to-moon legs, a duration of approximately 1.7 years, and a propellant usage of 125 kg), while offering greater observational performances and extending the mission to include Rhea, which was not addressed in previous studies. The four inter-moon transfers presented in this study require a total $\Delta V$ of 4356~m/s over approximately 3.2~years, corresponding to a propellant consumption of 231.5~kg. To highlight the relevance of the results obtained, reference transfers were computed using Hohmann impulsive maneuvers (with a specific impulse of 300~s) between coplanar moons, as well as planar 2BP circle-to-circle transfers assuming low-thrust applied tangentially at all times. The corresponding values are reported in Table \ref{tab:comparison}. The comparison clearly shows the advantage of exploiting the HIMs to exit and enter the SOI of each terminal moon. The HIM-based legs consistently reduce the $\Delta V$ with respect to the planar and unperturbed low-thrust reference solution, thus demonstrating the efficiency of the proposed strategy.

\begin{table*}[h!]
\centering
\small
\caption{Reference unperturbed planar transfers: LT circle-to-circle trajectories and Hohmann impulsive maneuvers, computed for an initial mass of 1000 kg.}
\label{tab:comparison}
\renewcommand{\arraystretch}{1}
\begin{tabular}{l|ccc|ccc}
\hline
{Transfer} &
\multicolumn{3}{c|}{{LT circle-to-circle}} &
\multicolumn{3}{c}{{Hohmann}} \\ 
\cline{2-7}
 & $\Delta V$ & $\Delta t$ & $\Delta m$ &
   $\Delta V$ & $\Delta t$ & $\Delta m$ \\
 & [m/s] & [day] & [kg] & [m/s] & [day] & [kg] \\
\hline
Rh-Di & 1542  & 472 & 93.5 & 1531 & 1.79 & 405 \\
Di-Te & 1320 & 407 & 73.2  & 1315 & 1.15 & 214 \\
Te-En & 1278 & 395 & 65.2 & 1274 & 0.81 & 134\\
En-Mi & 1675 & 510 & 77.8 & 1668 & 0.57 &  107 \\
\textbf{Total} & \textbf{5815} & \textbf{1784} & \textbf{309.7} & \textbf{5788} & \textbf{4.32} & \textbf{860} \\
\hline
\end{tabular}
\end{table*}

The tour architecture developed in this work differs fundamentally from those previously implemented or proposed in the literature. Cassini operated through a sequence of targeted flybys of the Saturnian moons using a chemical-propulsion system while traveling on highly elliptical Saturn-centered orbits. During its 13-year mission, Cassini completed 166 targeted flybys, 117 of which at Titan and 49 at the ILMs (23 of them at Enceladus) \cite{bellerose2019cassini}. Although the tour was designed to maximize scientific return, only a small fraction of each orbit was spent in the proximity of a moon, and no orbital capture was ever planned. During the orbital tour, approximately 700 kg of propellant were used, corresponding to a total $\Delta V$ of roughly 800 m/s, with an $I_{sp}$ of 323 s. Similarly, the strategy proposed by Strange et al.~\cite{strange2009leveraging} employs a propulsion system equivalent to that of Cassini and achieves a $\Delta V$ of 960 m/s with approximately 1000~kg of propellant. This tour enables 14 flybys at Rhea, 8 at Dione, 11 at Tethys, and 9 at Enceladus, but none at Mimas, over a total duration of 2.5 years. Moreover, based on previous work \cite{fantino2023end} on the interplanetary transfer, an unpowered capture at Titan can be achieved. In contrast, both the Cassini mission and the tour in \cite{strange2009leveraging} require a Saturn Orbit Insertion maneuver, with $\Delta V$ of about 600 m/s and 750 m/s, respectively.

In contrast, the present study adopts a completely different mission strategy. The tour is divided into five observational blocks, one for each ILM, starting at Rhea and ending at Mimas. Within each block, the spacecraft travels in a low-energy orbit for an arbitrary flight time, enabling prolonged and continuous observations with full-surface coverage. Relative velocities are of the order of 10–100 m/s, in contrast with velocities of the order of 1 km/s of Cassini and VILT-type tours, thus facilitating detailed scientific exploration. 
These weakly captured trajectories provide high inclinations, broad surface coverage, repeated close approaches, and negligible propellant cost, demonstrating the effectiveness of manifold-based dynamics for sustained lunar observation. After each science phase, the spacecraft transitions to the next moon through a low-thrust spiraling trajectory, resulting in long inter-moon transfer times. These legs are propagated in an ephemeris-based dynamics, an element absent from the literature on multi-moon tours. The $J_2$-perturbed CR3BP used for the intra-SOI phases naturally includes the dominant perturbation (oblateness of Saturn), which outweighs the cumulative third-body effects from the other moons and the Sun, whose influence remains below the adopted thresholds (1~km and 1~m/s). The result is a coherent and realistic representation of the Saturnian dynamical environment, allowing a weak-capture and low-thrust tour architecture not previously available in the literature. The remaining step toward a complete high-fidelity and implementable mission design is the inclusion of the orbital eccentricity and inclination of each moon in the intra-SOI dynamics, which will allow a complete unification of the dynamical models and enable a refined end-to-end mission scenario.

Another fundamental difference with respect to both Cassini and the existing literature is the use of electric propulsion. This choice provides high fuel efficiency and is selected to comply with strict power constraints at the distance of Saturn from the Sun. Although the limited thrust available inherently increases the inter-moon transfer times, the overall mission duration remains comparable to that of the primary and two extended Cassini missions. At the same time, the propellant requirements are significantly reduced. Consequently, the proposed architecture enables a complete multi-moon exploration at modest cost, in contrast to the higher fuel demands associated with chemical-propulsion and flyby-based mission concepts.

\section{Conclusions}\label{concl}
This work describes a mission design strategy for the exploration of the five Inner Large Moons (ILMs) of Saturn (Rhea, Dione, Tethys, Enceladus, and Mimas), combining natural dynamics with low-thrust propulsion. Motivated by the growing scientific interest in these bodies, the proposed concept balances the need for extended mission duration with substantial gains in propellant efficiency and scientific return.

In the present study, the exploration tour begins at Rhea and ends at Mimas. The interplanetary transfer from Earth to Saturn has been analyzed in~\cite{fantino2023end}, and is characterized by a total transfer duration of 12.34 years, from 21~January~2028 to 24~May~2040, and includes Saturn orbit insertion via a sequence of unpowered gravity assists at Titan. This interplanetary phase can be linked to the proposed three-dimensional tour, placing the current work within an end-to-end framework that extends from Earth departure to observations around Mimas.

The study at hand analyzes transfers between halo orbits around the libration points of the Saturn–moon-spacecraft three-body systems. The periodic nature of the departure and arrival orbits makes them effective gateways between consecutive science orbits and moon-to-moon transfers. Homoclinic and heteroclinic connections are obtained by leveraging the hyperbolic invariant manifolds (HIMs) associated with the periodic orbits, enabling a long-term weak capture around each target.  Connections between halo orbits of different moons are established with low-thrust arcs that bridge the distance and energy gap between the HIMs of consecutive moons.
This strategy demonstrates a strong scientific potential for the observation of the ILMs. Weakly captured three-dimensional orbits offer high inclinations, allowing polar regions observations, as well as repeated close approaches at negligible cost. Notably, visibility windows over the poles extend for several hours in Type B, C, and D heteroclinics, with the exact duration depending on the altitude. In contrast, homoclinic connections, due to their limited out-of-plane components, are shown to be less suitable for polar-region exploration. The dynamics within each sphere of influence (SOI) is modeled in the $J_2$–perturbed Circular Restricted Three-Body Problem (CR3BP), while the inter-SOI segments are propagated in a full $N$-body framework that also includes the oblateness of Saturn. The bodies considered in the $N$-body model are selected through a dedicated perturbation analysis, and their ephemeris extracted from SPICE, ensuring dynamical fidelity while maintaining computational efficiency.

In terms of performance, the four inter-moon legs require a total $\Delta V$ of 4.4 km/s over approximately 3.2~years, corresponding to a propellant consumption of 231.5~kg for an initial mass of 1000~kg. The duration of the tour is further extended through science-orbit dwelling around each moon, which can be swept for arbitrary lengths at marginal propellant cost. The results clearly demonstrate the benefits of integrating CR3BP-based dynamical structures with low-thrust propulsion, providing a powerful framework for high-fidelity and fuel-efficient transfer design for multi-moon environments. The proposed design strategy is compatible with current propulsion and power technologies and supports the scientific and operational objectives of a future ILMs exploration mission. The methodology presented is not limited to the Saturnian system and may be extended to multi-moon systems of other giant planets, which represent compelling targets for future scientific exploration.

Future work will focus on incorporating the orbital eccentricity and inclination of each ILM into the intra-SOI dynamics, a step toward achieving a high-fidelity design across all phases of the mission. 
Given the advantages provided by the science orbits identified in this study in terms of coverage and cost, the next refinement is to transition them within a more accurate dynamical framework that admits the same manifold-based structures. For this reason, a future study involves the transition to a $J_2$-perturbed Elliptic Restricted Three-Body Problem (ER3BP), which preserves the dynamical features exploited in the  current work in a more realistic model.
Once the ER3BP-based orbits are generated, the effects of the inclinations of the lunar orbits and additional third-body perturbations will be compensated through corrective maneuvers, which would allow the trajectories to follow their nominal manifold-driven ones. These adjustments are expected to require minimal propellant consumption, enabling the weak-capture and low-thrust architecture developed in this work to evolve into a fully realistic design suitable for future mission implementation.

\section*{Acknowledgment}
This research was supported by Polar Research Center (PRC), Khalifa University of Science and Technology (KU-PRC).
The authors thank Dr. Roberto Flores for the insightful discussions regarding the perturbation analysis.


\end{document}